\documentclass[a4paper,12pt]{article} 
 \usepackage[a4paper,top=0.75in,bottom=1.15in,left=0.75in,right=0.75in,marginparwidth=0.75cm]{geometry}
\usepackage{amsmath}
\usepackage{amssymb}
\usepackage{mathtools} 
\usepackage{enumerate}
\usepackage{graphicx}
\usepackage{float}
\usepackage{cite}
\usepackage{color}
\usepackage{url}
\usepackage{subcaption}
\usepackage{caption}
\usepackage{hyperref}
\usepackage{float}
\usepackage{siunitx}
\usepackage{booktabs}
\usepackage[english]{babel}
\usepackage[nottoc]{tocbibind}
\usepackage{lipsum}
\usepackage{booktabs}
\usepackage{siunitx}
\usepackage{makecell}
 
\usepackage{changepage}
\usepackage{orcidlink}
\usepackage{hyperref}
\definecolor{DarkBlue}{rgb}{0.00,0.00,0.50}
\hypersetup{colorlinks=true,
            linkcolor=blue,
            filecolor=magenta,
            urlcolor=cyan,
            citecolor=DarkBlue,
            pdftitle={MSc UoR Computer Science Report},
            }
\usepackage{url}
\usepackage{caption}
\usepackage{array}
\usepackage{rotating}
\usepackage{algorithm,algorithmic}
\usepackage{amsthm}
\newtheorem{definition}{Definition}

\newtheorem{defn}{Definition}

\newtheorem{prb}{Test}
\let\oldprb\prb
\renewcommand{\prb}{\oldprb\normalfont}
\newtheorem{lemma}{Lemma}
\let\oldlemma\lemma
\renewcommand{\lemma}{\oldlemma\normalfont}
\newtheorem{pol}{Polynomial}
\let\oldpol\pol
\renewcommand{\pol}{\oldpol\normalfont}
\usepackage{authblk}
\usepackage{blindtext}
  \begin{document}
\title{Modified approach for linear and non-linear IBVPs with fractional dynamics}
\date{}
 \author{ Qasim Khan$^{*}$}
\author{   Anthony Suen}
\affil{\normalsize{ Department of Mathematics and Information Technology,\\  The Education University of Hong Kong,\\ 10 Lo Ping Road, Tai Po, N.T,\\ Hong Kong} }
\affil{\textit {qasimkhan@s.eduhk.hk  \ \
acksuen@eduhk.hk}}
\date{\today}
\maketitle
 	
\begin{abstract}
Analytical and numerical techniques have been developed for solving fractional partial differential equations (FPDEs) and their systems with initial conditions. However, it is much more challenging to develop analytical or numerical techniques for FPDEs with boundary conditions, although some methods do exist to address such problems. In this paper, a modified technique based on the Adomian decomposition method with Laplace transformation is presented, which effectively treats initial-boundary value problems. The non-linear term has been controlled by Daftardar-Jafari polynomials.  Our proposed technique is applied to several initial and boundary value problems and the obtained results are presented through graphs. The differing behavior of the solutions for the suggested problems is observed by using various fractional orders. It is found that our proposed technique has a high rate of convergence towards the exact solutions of the problems. Moreover, while implementing this modification, higher accuracy is achieved with a small number of calculations, which is the main novelty of the proposed technique. The present method requires a new approximate solution in each iteration that adds further accuracy to the solution. It demonstrates that our suggested technique can be used effectively to solve initial-boundary value problems of FPDEs.
\end{abstract}
\textbf{Keywords:} {Initial and boundary value problems}; {  Decomposition method}; Laplace transformation; {Caputo derivative}; Daftardar-Jafari polynomials; Modified Laplace decomposition method.

\tableofcontents
\section{Introduction}
In recent decades, fractional partial differential equations (FPDEs) have gained significant attention in research due to their robust modeling capabilities and wide-ranging applicability across various scientific and engineering fields. As a generalization of classical partial differential equations (PDEs), FPDEs allow researchers to effectively capture complex behaviors and phenomena that standard models may not fully address. This versatility is evident in diverse applications, including electrochemistry, where they model charge transport processes, and turbulence flow, where they describe intricate fluid dynamics. FPDEs are also employed in material sciences, chaotic dynamics, diffusion processes, electromagnetics, ecology, epidemiology, economic systems, quantum mechanics, and magnetohydrodynamics, among others. Furthermore, their significance extends to medical imaging, traffic flow modeling, and the electrical signaling of nerves, highlighting their potential to tackle real-world challenges across multiple domains. The growing interest in FPDEs reflects their role as a powerful tool for researchers seeking to develop more sophisticated and effective models of complex systems. For more details, see the references cited in \cite{10,11,12,13,14,15,16}.

Initial and boundary value problems (IBVPs) for FPDEs are particularly important as they arise in practical applications where the state of a system is influenced by both initial conditions and constraints at the boundaries of the domain. The solutions to these problems provide critical insights into the behavior of dynamic systems over time and space. The mathematical treatment of IBVPs for FPDEs is often more complex than their integer-order counterparts, due to the non-local nature of fractional derivatives and the need for appropriate boundary conditions that reflect the physical reality of the modeled phenomena. Additionally, the interplay between various types of boundary conditions, the behavior of solutions in different regions of the domain, and the presence of nonlinearity further complicate the analysis and computation of PDEs, making it an active area of research in applied mathematics and engineering.

The study of initial and boundary value FPDEs encompasses a range of methodologies, including analytical techniques, numerical simulations, and qualitative analyses. Researchers have devised a variety of approaches to address the complexities associated with FPDEs, for example fractional finite difference methods, spectral methods, and variational techniques, which are employed to obtain approximate solutions and to gain insights into the properties of these equations. Notable techniques in this field include the wavelet optimization method \cite{17}, Bäcklund transformation technique \cite{18}, moving finite element method \cite{19}, finite difference method \cite{20}, first integral method \cite{21}, modified exp-function method \cite{22}, spectral–Galerkin method \cite{23}, stable multi-domain spectral penalty method \cite{24}, Haar wavelet Picard method \cite{25}, neural network method \cite{26}, spectral Petrov‐Galerkin method \cite{27}, high-order finite element method \cite{28}, modified trial equation method \cite{29}, residual power series method \cite{H3}, block pulse operational matrix method \cite{30}, Laplace inversion technique \cite{MLADM45}, weighted finite difference methods \cite{32}, wavelet operational method \cite{33}, Yang-Laplace decomposition method \cite{34}, Aboodh transform  decomposition method \cite{aboodh1,aboodh2}, new approximate analytical method \cite{H2}, hybrid natural transform homotopy perturbation method \cite{35}, meshless technique \cite{36}, invariant subspace method \cite{37}, and sumudu decomposition method \cite{38}.

 The Adomian Decomposition Method (ADM), created by George Adomian in the 1980s, is recognized as a highly effective technique for addressing nonlinear functional equations. An example of its application can be found in the work of Li and Pang \cite{48}, who applied this method to solve nonlinear systems. S.S. Ray in \cite{49} utilized ADM to solve the Bagley-Torvik equation, while J. Biazar in \cite{50} employed the same technique to tackle a system of integral-differential equations. In \cite{57}, A.M. Waz Waz solved time-dependent Emden–Fowler type equations and Qibo Mao \cite{58} studied the free vibration analysis of multiple-stepped beams. Different researchers have also applied modifications of ADM to solve FPDEs, including a new modification of ADM that was introduced by Elaf Jaafar Ali for solving integer-order IBVPs in \cite{60,59} and \cite{ali}.

 In this research article, we introduce a novel modification of decomposition methods that is integrated with the Laplace transformation to tackle both linear and nonlinear initial boundary value problems involving fractional orders. Our modified approach aims to improve the effectiveness and precision in solving complex fractional differential equations, which frequently arise in diverse fields such as physics, engineering, and applied mathematics. By combining the benefits of modified decomposition techniques and the analytical power of the Laplace transformation, we provide a robust framework for addressing these intricate mathematical challenges. In the literature, there are only very few methods to solve IBVPS. So in this connection, we try to extend the contribution of Elaf Jafaar Ali for solving fractional order FPDEs with IBVPs. By applying this technique to IBVPs, we compute a new approximate initial solution with a twice-recursive formula for each iteration, which can enhance the rate of convergence of the method. It is worth noting that if the FPDEs or their corresponding system have higher number of initial boundary conditions, then a high number of calculations are required to obtain the desired results. In contrast, the accuracy of the suggested technique is cleared from thing that very small number of calculations are needed for computation, and our results are illustrated through plots.

 The rest of the paper is organized as follows: Section \ref{sec2} demonstrates the proposed research methodology, Section \ref{sec3} presents some experiments on the proposed method, and the obtained results are discussed in Section \ref{sec4}.  In the last section of the paper \ref{conclusion}, we summarize our findings and suggest several avenues for future research. This discussion highlights areas where further exploration could enhance understanding or application of our proposed methods.

\section{Research  Methodology}\label{sec2}
To understand the basics of our proposed research methodology, we split Section \ref{sec2} into four subsections. Subsection \ref{subsec2.1} provides an overview of the traditional Adomian Decomposition Method. Then, we recall the modification of ADM in Subsection \ref{subsec2.2}. In the same manner, we present the classical methodology of the Laplace Adomian Decomposition Method and modify it for nonlinear initial and boundary value problems in Subsections \ref{subsec2.3} and \ref{subsec2.4}, respectively.
\subsection{Adomian Decomposition Method} \label{subsec2.1}
The Adomian Decomposition Method (ADM) is a semi-analytical approach used to solve ordinary and partial nonlinear differential equations. This innovative method was developed by George Adomian from the 1970s to the 1990s \cite{adomian1988review,adomian2013solving}. ADM has become an effective tool for solving both linear and nonlinear fractional PDEs in recent years. To understand the fundamental concept behind this method, we consider a general equation of the form:
\begin{equation}\label{L2}
\mathrm{L}{{u}} + Q{{u}} + {{{N}}}{{u}} = {{h}},
\end{equation}
where $\mathrm{L}$ is a first-order linear operator that can be easily inverted, $Q$ represents a linear term, ${{{N}}}{{u}}$ is a nonlinear operator, and $g$ is the source term.\\
 To derive the solution, we apply the inverse ($\mathrm{L}^{-1}$) to each sides of Eq. \eqref{L2}:
\begin{equation}\label{F1}
{{u}} = {{{c}}} + \mathrm{L}^{-1}({{h}}) - \mathrm{L}^{-1}(Q{{u}}) - \mathrm{L}^{-1}({{{N}}}{{u}}),
\end{equation} 
In this context, the variable  ${{c}}$
 denotes the constant of integration. The solution derived from the Adomian Decomposition Method (ADM) is formulated in the form of a series expansion, which facilitates the approximation of the desired function. This series representation allows for the systematic computation of terms, each of which contributes to the overall solution. 
  \begin{equation}\nonumber
{{u}} = \sum_{n=0}^{\infty}{{u}_{n}}.
\end{equation}
 The nonlinear term ${{{N}}}{{u}}$ has been controlled by Adomian polynomials and represented by the series
 $A_n$\\
\begin{equation}\label{srs}
{{{N}}}{{u}} = \sum_{n=0}^{\infty} A_n.
\end{equation} 
The definition of  $A_n$ is provided below.
\begin{definition}\label{def AP}
The function $A_n$ as appeared in Eq. \eqref{srs} is known as the Adomian polynomial:
  \begin{equation}\nonumber
A_n = \frac{1}{n!} \frac{d^{n}}{d{{\lambda}}^{n}} {{{N}}}\left(\sum_{n=0}^{\infty} ({{\lambda}}^{k}{{u}}_{n})\right), \quad n=0,1,\dots
\end{equation}
\end{definition}

By applying Adomian polynomial, we provide the iteration formula for the solution $u$ as given in Eq. \eqref{F1}:
\begin{equation}\nonumber
\begin{split} 
{{u}}_{0} = {{{c}}} + \mathrm{L}^{-1}({{h}}), \quad \textnormal{initial iteration}.\\
 {{u}}_{n+1} = \mathrm{L}^{-1}(Q{{u}}_{n}) - \mathrm{L}^{-1}(A_n), \quad n \geq 0.
\end{split}
\end{equation}
The recursive relations of the ADM rely solely on the initial condition and the source term. 
\subsection{Modified ADM for Linear IBVPs}\label{subsec2.2}
The basic idea discussed in Subsection \ref{subsec2.1} can been applied to initial value problems for partial differential equations. In recent years, a number of researchers have successfully extendeded this idea to linear initial and boundary value problems with fractional order (see \cite{M01,M1,M2,M3,M4,M5}). To illustrate the method, we consider the following one-dimensional linear equation
\begin{equation}\label{eq1}
\frac{\partial^{\alpha}}{\partial t^{\alpha}} u(x, t) + u(x, t) = \frac{\partial^{2} }{\partial x^2} u(x, t) + h(x, t), \quad 0 \leq x \leq 2, \quad 0 < \alpha \leq 1, \quad t \geq 0,
\end{equation}
which is equipped with the initial and boundary conditions (IBCs) as follows:
\begin{equation}\nonumber
\begin{split} 
u(x, 0) &= {{{f}}} (x), \quad 0 \leq x \leq 2,\\
u(0, t) &= {{{g}}}_{0}(t), \quad \\
u(2, t) &= {{{g}}}_{1}(t), \quad t > 0.
\end{split}
\end{equation}
In Eq. \eqref{eq1}, the source term is represented by $h(x, t)$. In operator form, Eq. \eqref{eq1} can be further rewritten as
\begin{equation}\label{EQ}
\mathrm{L} u = \frac{\partial^{2} }{\partial x^2}u(x, t) - u(x, t) + h(x, t),
\end{equation}
where the operator $\mathrm{L}$ is defined as
\begin{equation}\nonumber
\mathrm{L} = \frac{\partial^{\alpha}}{\partial t^{\alpha}}.
\end{equation}
We define the inverse operator $\mathrm{L}^{-1}$ by
\begin{equation}\nonumber
\mathrm{L}^{-1}(\cdot) = I^{\alpha}(\cdot) d t,
\end{equation}
where the operator $I^{\alpha}$ is a fractional integral. The following lemma gives some properties of $I^{\alpha}$ and the proof can be found in reference \cite{51}.
\begin{lemma}\label{lemma1}
   For $ j-1 < \alpha \leq j$ with $j \in \mathbb{N}$ and $ h : \mathbb{R}^+\to\mathbb{C}$, the operator $I^{\alpha}$ satisfies the following identities
   \begin{equation}
\begin{cases}
I^{\alpha} I^{b} h(x) &= I^{\alpha + b} h(x), \quad  b, \alpha \geq 0,\\
I^{\alpha} x^{\lambda} &= \frac{\Gamma(\lambda + 1)}{\Gamma(\alpha + \lambda + 1)} x^{\alpha + \lambda}, \quad \alpha > 0,\,\, \lambda > -1,\,\, x > 0,\\
I^{\alpha} D^{\alpha} h(x) &= h(x) - \sum_{k=0}^{j-1} h^{(k)}(0^{+}) \frac{x^{k}}{k!}.
\end{cases}
\end{equation}
\end{lemma}
Applying the operator $\mathrm{L}^{-1}$ to Eq. \eqref{EQ}, we obtain
\begin{equation}\nonumber
u(x, t) = u(x, 0)   + \mathrm{L}^{-1}\left(\frac{\partial^{2}  }{\partial x^{2}}u(x, t) -  u(x, t) + h(x, t)\right).
\end{equation}
The initial approximate solution becomes
\begin{equation}\nonumber
u_{0}(x, t) = u(x, 0) +  \mathrm{L}^{-1}(h(x, t)),
\end{equation}
and the iterative formula is expressed as follows:
\begin{equation}\nonumber
u_{n+1}(x, t) = \mathrm{L}^{-1}\left(\frac{\partial^{2} }{\partial x^{2}}u_{n}^{*}(x, t) - u_{n}^{*}(x, t)\right).
\end{equation}
The new   approximate solution $u_{n}^{*}$ for $u(x,t)$ can be computed as
\begin{equation}\label{ali-operator} 
u_{n}^{*} = u_{n}(x, t) + (1 - x)\left[u(0,t) - u_{n}(0, t)\right] + x\left[u(2,t) - u_{n}(2, t)\right].
\end{equation}
The new solutions $u_{n}^{*}$ of Eq. \eqref{eq1} satisfy both the IBCs, as given below:
\begin{equation}\nonumber
\begin{split}
&\text{at } t = 0, \quad u_{n}^{*}(x, 0) = u_{n}(x, 0),\\
&\quad x = 0, \quad u_{n}^{*}(0, t) = u_{0}(t),\\
&\quad x = 2, \quad u_{n}^{*}(2, t) = u_{1}(t).
\end{split}
\end{equation} 
The more general case of $u_{n}^{*}\backsim$ \ref{ali-operator}  will be discussed in the modified Laplace decomposition method in subsection \ref{subsec2.4}.
\subsection{Classical LADM for  IVPs}\label{subsec2.3}
Before proceeding to achieved our main goal, we will now consider the following initial value generalized equation
\begin{equation}\label{eq11}
 {}^{C}_{}D^{\alpha}_{t} u({{{x}}}, {{{t}}}) + Qu({{{x}}}, {{{t}}}) + Nu({{{x}}}, {{{t}}}) = h({{{x}}}, {{{t}}}), \quad {{{x}}}, {{{t}}} \geq 0, \quad n - 1 < \alpha \leq n,
\end{equation}
equipped with the initial condition
\begin{equation*}
u({{{x}}}, 0) = f({{{x}}}).
\end{equation*}
The fractional derivative in Eq. \ref{eq11} is expressed in the Caputo sense, which will be given in Def. \ref{def2}. The linear and nonlinear terms are denoted by $ Q$   and $N$ respectively, and $h({{{x}}}, {{{t}}})$ is a given source term.
\begin{definition}\label{def2}
The Caputo fractional derivative is defined as a way to extend the idea of taking derivatives to non-integer orders. It uses an integral operator to allow for the differentiation of functions, which helps in studying systems with memory effects and past influences. This type of derivative is especially useful when dealing with initial conditions that are defined in a traditional manner, making it applicable in many areas of science and engineering.
\begin{equation}\label{cd}
  {}^{C}_{}D^{\alpha}_{t}{u}({{{x}}},{{{t}}}) = \frac{1}{{\Gamma}(n-\alpha)} \int_{0}^{t} (t - \tau)^{n - \alpha - 1} \frac{\partial^n}{\partial \tau^n} u(x,\tau) \, \mathrm{d}\tau
\end{equation}
where $n - 1 < \alpha \leq n$ and $n \in \mathbb{N}$. More discussion can be found in \cite{caputo1967linear}.
 \end{definition}
\begin{definition}
The Laplace transform of the Caputo fractional derivative can be expressed as follows (\ref{lcd}): it provides a way to convert the fractional derivative into a different form that is easier to work with, especially in solving differential equations. This transformation helps simplify calculations by turning complex functions into algebraic equations in the Laplace domain
\begin{equation}\label{lcd}
  \mathcal{L} \left( {}^{C}_{}D^{\alpha}_{t} u({{{x}}},{{{t}}}) \right) = s^{\alpha} u(x,s) - \sum_{k=0}^{n-1} s^{\alpha-k-1} \frac{\partial^{k}}{\partial t^{k}} u(x,0)
\end{equation}
where $u(x,s) = \mathcal{L}(u(x,t))$. If $n=1$, Eq. \eqref{lcd} can be simplified to:
\begin{equation}\label{lcd1}
  \mathcal{L} \left( {}^{C}_{}D^{\alpha}_{t} u({{{x}}},{{{t}}}) \right) = s^{\alpha} u(x,s) - s^{\alpha - 1} u(x,0)
\end{equation}
where $0 < \alpha \leq 1$.
 \end{definition}
Applying Laplace transform $\mathcal{L}$ on both sides of Eq. \eqref{eq11}, we obtain
\begin{equation*}
s^\alpha \mathcal{L}\{u({{{x}}}, {{{t}}})\} - s^{\alpha - 1} u({{{x}}}, 0) = \mathcal{L}\{h({{{x}}}, {{{t}}})\} + \mathcal{L}\{Qu({{{x}}}, {{{t}}}) + Nu({{{x}}}, {{{t}}})\},
\end{equation*}
Re-arranging for $\mathcal{L}\{u({{{x}}}, {{{t}}})\}$, we obtained 
\begin{equation}\label{p3}
\mathcal{L}\{u({{{x}}}, {{{t}}})\} = \frac{u({{{x}}}, 0)}{s} +\frac{1}{s^\alpha}  \mathcal{L}\{Qu({{{x}}}, {{{t}}}) + Nu({{{x}}}, {{{t}}}) + h({{{x}}}, {{{t}}})\}.
\end{equation}
We let $u(x,t)$ to be the decomposed  solution, which is given by
\begin{equation}\label{p4}
u({{{x}}}, {{{t}}}) = \sum_{j=0}^{\infty} u_j({{{x}}}, {{{t}}}).
\end{equation}
The nonlinear term in Eq. \eqref{eq11} can be expressed as
\begin{equation}\label{p5}
Nu({{{x}}}, {{{t}}}) = \sum_{j=0}^{\infty} A_j,
\end{equation}
where $A_j$ is the Adomian polynomials defined in Def. \eqref{def AP}.

By substituting Eq. (\ref{p4}) and Eq. (\ref{p5}) into Eq. (\ref{p3}), we obtain:
\begin{equation}\label{eq32}
\mathcal{L} \left\{  \sum_{j=0}^{\infty} u({{{x}}}, {{{t}}}) \right\} = \frac{f({{{x}}})}{s} +\left( \frac{1}{s^\alpha} \mathcal{L} \left\{\sum_{j=0}^{\infty} Q u_j({{{x}}}, {{{t}}})\right\} + \sum_{j=0}^{\infty} A_j + h({{{x}}}, {{{t}}})\right).
\end{equation}
Applying the decomposition method, we get
\begin{equation}\label{eq33}
\mathcal{L}\{u_0({{{x}}}, {{{t}}})\} =  \frac{f({{{x}}})}{s} +\frac{1}{s^\alpha}\mathcal{L} \{h({{{x}}}, {{{t}}})\}.
\end{equation}
and
\begin{equation}\label{eq34}
\mathcal{L}\{u_{j+1}({{{x}}}, {{{t}}})\} =   \mathcal{L}^{-1}\frac{1}{s^\alpha}\mathcal{L}\{Qu_j({{{x}}}, {{{t}}}) + A_j\}, \quad j \geq 1.
\end{equation}
Using the inverse transform on equations (\ref{eq33}) and (\ref{eq34}), we get
\begin{equation*}
u_0({{{x}}}, {{{t}}}) = f({{{x}}})+\mathcal{L}^{-1}\left\{\frac{1}{s^\alpha}\mathcal{L} \{h({{{x}}}, {{{t}}})\} \right\},
\end{equation*}
and
\begin{equation*}
u_{j+1}({{{x}}}, {{{t}}}) =  \mathcal{L}^{-1}\left\{ \frac{1}{s^\alpha} \mathcal{L}\left\{ Q u_j({{{x}}}, {{{t}}})  + A_j\right\} \right\}.
\end{equation*}
The recursive relation is based on the first iteration, which relies solely on the initial condition with the source term. The method has issues when boundary conditions are involved. There is no doubt that classical LADM has many applications in solving linear and nonlinear initial value problems (refer to refs. \cite{LADM1,LADM3,ongun2011laplace,kiymaz2009algorithm} and the references cited therein).
\subsection{Modified LADM for non-linear IBVPs} \label{subsec2.4}
To better understand the key concepts of the modified LADM, let us examine the following generalized equation related to initial and boundary value problems
\begin{equation}\label{eq111}
 {}^{C}_{}D^{\alpha}_{t} u({{{x}}}, {{{t}}}) + Qu({{{x}}}, {{{t}}}) + Nu({{{x}}}, {{{t}}}) = h({{{x}}}, {{{t}}}),\quad l \leq {{{x}}} \leq L, \quad  {{{t}}} \geq 0, \quad 0 < \alpha \leq 1
\end{equation}
with the initial and boundary conditions (IBCs) as follows: 
\begin{equation*}
\begin{split}
&u({{{x}}}, 0) = f({{{x}}}), \\
&u({{{l}}}, t) = g_{0}, \\
&u({{{L}}}, t)  = g_{1}.
\end{split}
\end{equation*}
The fractional derivative in Eq. \eqref{eq111} is expressed in the Caputo sense as defined in Def. \ref{def2}. The linear and nonlinear terms are denoted by $ Q$ and $N$ respectively, and $h({{{x}}}, {{{t}}})$ is the source term. We also denote $ l $ as the lower boundary limit and $ L $ as the upper boundary limit in Eq. \ref{eq111}.

It is clear from the classical LADM that the first iteration takes the following form:
\begin{equation*}
u_0({{{x}}}, {{{t}}}) = f({{{x}}})+\mathcal{L}^{-1}\left\{\frac{1}{s^\alpha}\mathcal{L} \{h({{{x}}}, {{{t}}})\} \right\}.
\end{equation*}
The new approximate solution, denoted by $u^{*}(x,t)$, was proposed by Elaf Jaafar in \cite{59,ali} to solve  linear initial boundary value problems (IBVPs) and is expressed as follows:
\begin{equation}\label{qwas}
u_{n}^{*}(x,t) = u_{n}(x, t) + (1 - x)\left[g_{0} - u_{n}(l, t)\right] + x\left[g_{1} - u_{n}(L, t)\right].\tag{*}
\end{equation}
It is clear that the modified approximate solution fulfills both the initial and boundary conditions, as demonstrated below:
\begin{equation}\nonumber
\begin{cases}
&if \quad t = 0, \quad u_{n}^{*}(x, 0) = u_{n}(x, 0),\\
&if \quad x = l, \quad u_{n}^{*}(l, t) = u_{n}(l, t),\\
&if \quad x = L, \quad u_{n}^{*}(L, t) = u_{n}(L, t).
\end{cases}
\end{equation}
The authors in \cite{lagaris1998artificial} employed the operator defined in eq (\ref{qwas}) within their trial solutions, utilizing artificial neural networks to effectively solve ordinary and partial differential equations, as reported in 1997. Due to the effectiveness of this operator, researchers are now applying it in machine learning and artificial neural networks. The most up-to-date reference can be found in \cite{3ali2024adaptive}, and we recommend that interested readers consult these sources as well as the references cited therein.

Finally, we propose a modified recursive relation which can be expressed as follows
\begin{equation*}
u_{n+1}({{{x}}}, {{{t}}}) =  \mathcal{L}^{-1}\left\{ \frac{1}{s^\alpha} \mathcal{L}\left\{ Q u^{*}_n({{{x}}}, {{{t}}})  + B^{*}_n\right\} \right\},  \forall n \in \mathbb{N} .
\end{equation*}
where $B^{*}_n=N \left(\sum_{{n}=0}^\infty{{{u}_{n}^{*}}(x,{{{t}}})}\right) $ that will be given in the following definition.
\begin{definition}\label{D-J}
\textnormal{ The function $B^*_n$ is known as the Jafari polynomials which is given by\\
\begin{equation}\label{s4}
B^*_n=N \left(\sum_{{n}=0}^\infty{{{u}_{n}^{*}}(x,{{{t}}})}\right)=N({{{u}}_0^{*}}(x,{{{t}}})+\sum_{{n}=0}^\infty \bigg[{N} \bigg(\sum_{i=0}^{{n}} {{{u}}_i^{*}}(x,{{{t}}})\bigg)-{N} \bigg(\sum_{i=0}^{{n}-1 }{{u}}_i^{*} (x,{{{t}}})\bigg)\bigg].
\end{equation}
We refer the interested readers to \cite{cost105,cost1} for more details.}
\end{definition}

The key steps for our Modified Laplace
 Decomposition Method (MLDM) can be summarized as follows:
\renewcommand{\thealgorithm}{}
 \begin{algorithm*}\nonumber
    \caption{Modified Laplace Decomposition Technique within Jafari Polynomials}
    \begin{algorithmic}[1]
        \STATE Define the fractional derivative in the Caputo sense with its Laplace transformation.
        \STATE Use the decomposition procedure to simplify the algebraic equation.
        \STATE Use the initial guess and source term for the first iteration, $u_{0}(x,t)$.
        \STATE Obtain a new initial approximation by setting $n := 0$ and assigning $g_{0}$ and $g_{1}$ to the code: i.e.,
        $
        u_{0}^{*}(x,t) := u_{0}(x, t) + (1 - x)\left[g_{0} - u_{0}(l, t)\right] + x\left[g_{1} - u_{0}(L, t)\right].
         $
        \STATE Control the non-linear term using Jafari polynomials.
        \STATE Utilize $u_{0}^{*}(x,t)$ to obtain the second iteration,\\
         i.e., for $n = 0$,
        $
        u_{1}(x,t) := \mathcal{L}^{-1}\left\{ \frac{1}{s^\alpha} \mathcal{L}\left\{ Q u_{0}^{*}(x,t) + B_{0}^{*} \right\} \right\}.
        $
        \STATE Repeat for all $n \in 1, 2, 3,\cdots$.
        \STATE Use the inverse Laplace transform to convert the expression from the $s$-domain to the time domain.
        \STATE Sum the first few terms of the series to obtain an approximate solution.
        \STATE Compare the obtained approximate solution with the classical Laplace Adomian Decomposition Method (LADM) to test the accuracy.
    \end{algorithmic}
    \label{al1}	
\end{algorithm*}
\section{ Numerical Experiments }\label{sec3}
In this section, we demonstrate the effectiveness and reliability of the Modified Laplace Adomian Decomposition Method (MLDM) by applying it to several experimental tests. The problems we are investigating has been previously examined by numerous researchers, as referenced in works such as \cite{MLADM45,moaddy2011non,odibat2009variational,1yavuz2018numerical,12chen2010fractional,3ali2024adaptive} and some related literature of initial and boundary value fractional order differential equations. In this study, we apply the MLDM to compare its outcomes with those obtained from the classical LADM. While most of our approximate solutions are based on a two-term approximation, we also extend the series solution to four terms in certain cases to enhance accuracy.
\begin{prb}\label{problem1}
We consider the following one-dimensional fractional PDEs:
\begin{equation}\label{EQ1}
\frac{\partial^{\alpha} }{\partial t^{\alpha}}u(x, t) + u(x, t) = \frac{\partial^{2} }{\partial x^{2}}u(x, t) + h(x, t), \quad 0 \leq x \leq 2, \quad 0 < \alpha \leq 1, \quad t \geq 0,
\end{equation}
with IBCs as follows:
\begin{equation}\nonumber
\begin{split}
u(x, 0) &= 0, \\
u(0, t) &= u(2, t) = 0,\\
\end{split}
\end{equation}
where the source term  $h(x,t)$ is given by
\begin{equation*}
h(x, t) = \frac{2}{\Gamma(3 - \alpha)} x (2 - x) + 2t^{2}.
\end{equation*}
The analytical solution is given by
\begin{equation}\nonumber
u(x, t) = t^{2} x (2 - x).
\end{equation}
\begin{figure}[H]
    \centering
    \includegraphics[width=5.4cm]{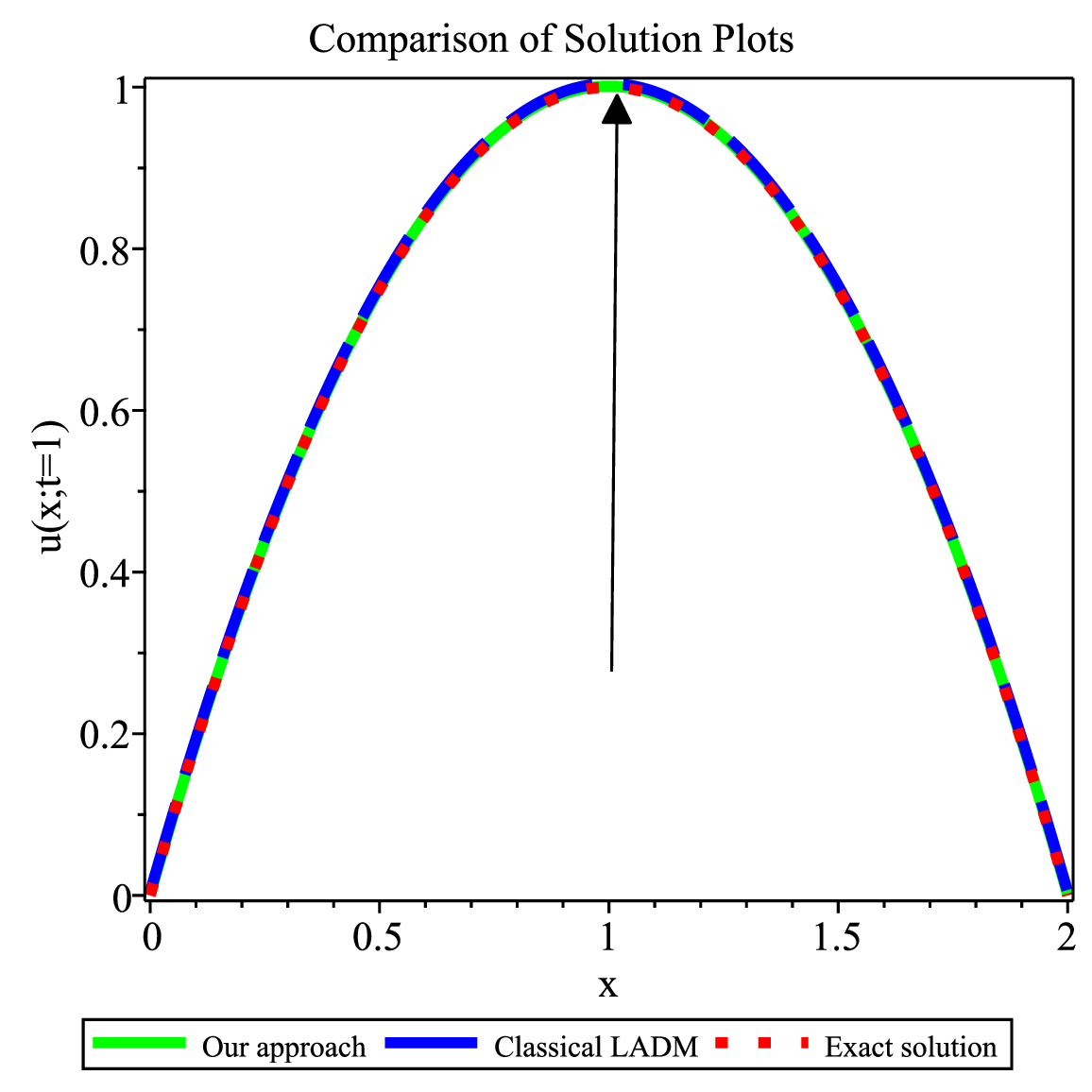}%
    \llap{\raisebox{0.95cm}{%
      \includegraphics[width=2.50cm]{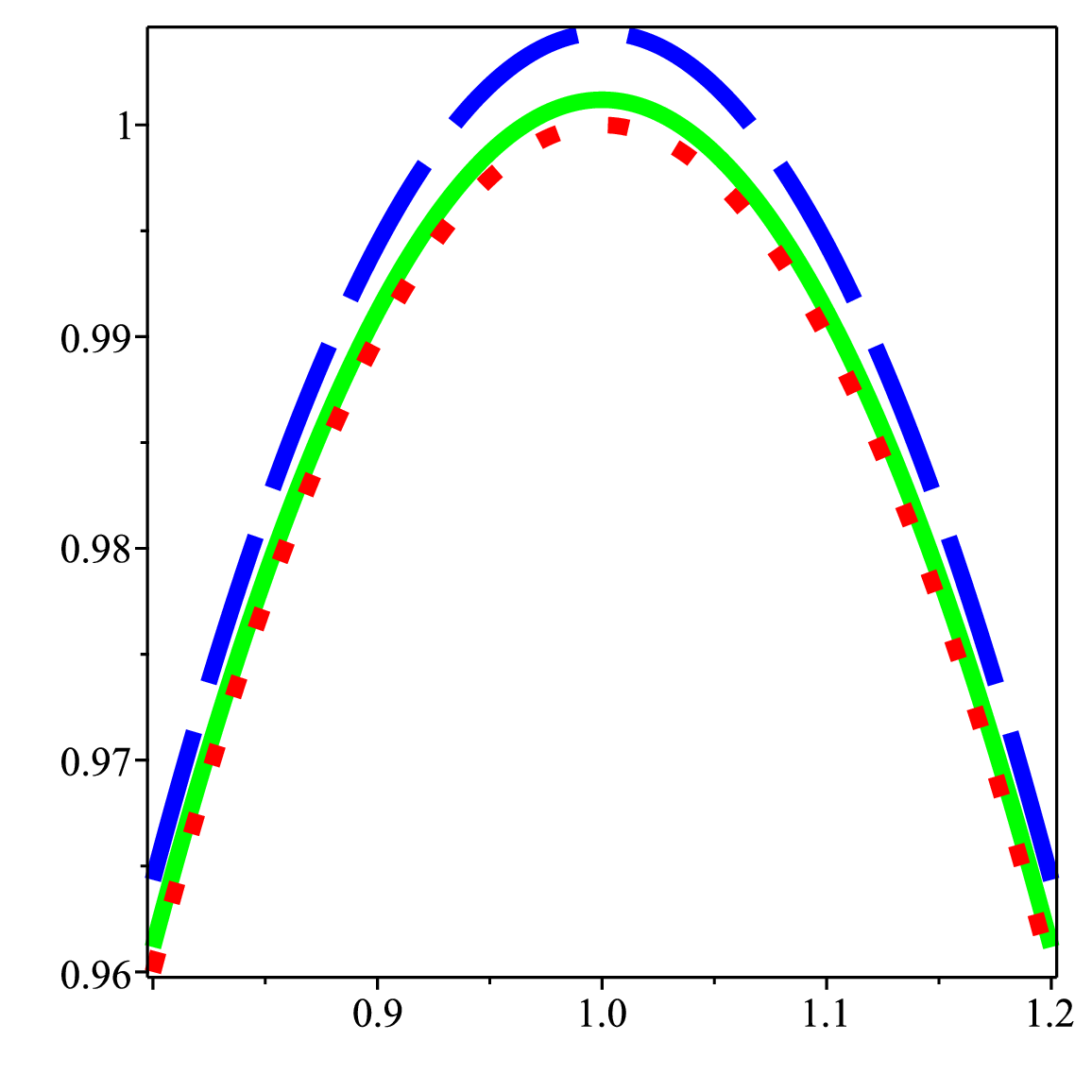}\hspace{2.90em}%
    }}\hspace{0.1em}
     \includegraphics[width=5.4cm]{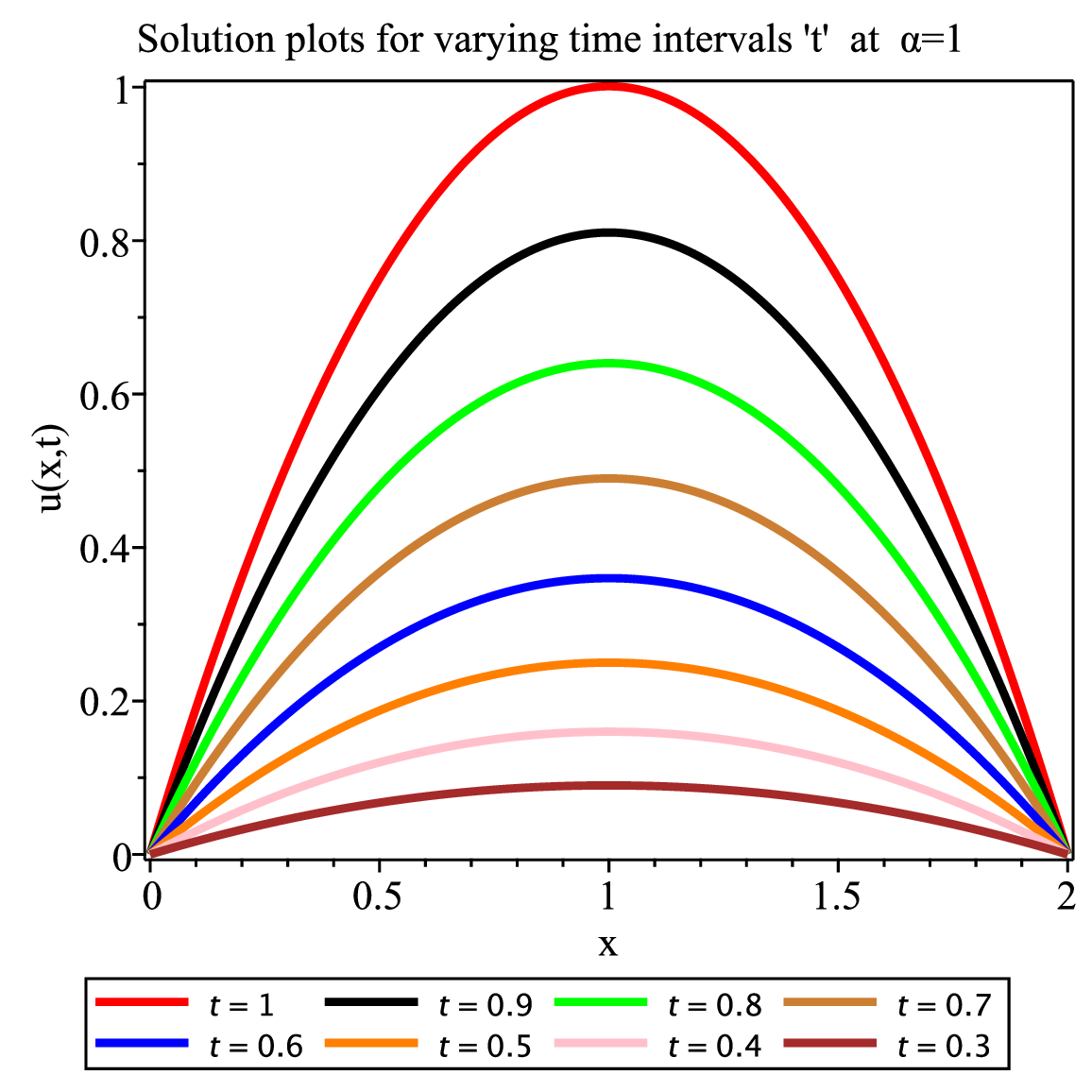}\hspace{0.1em}%
     \includegraphics[width=5.4cm]{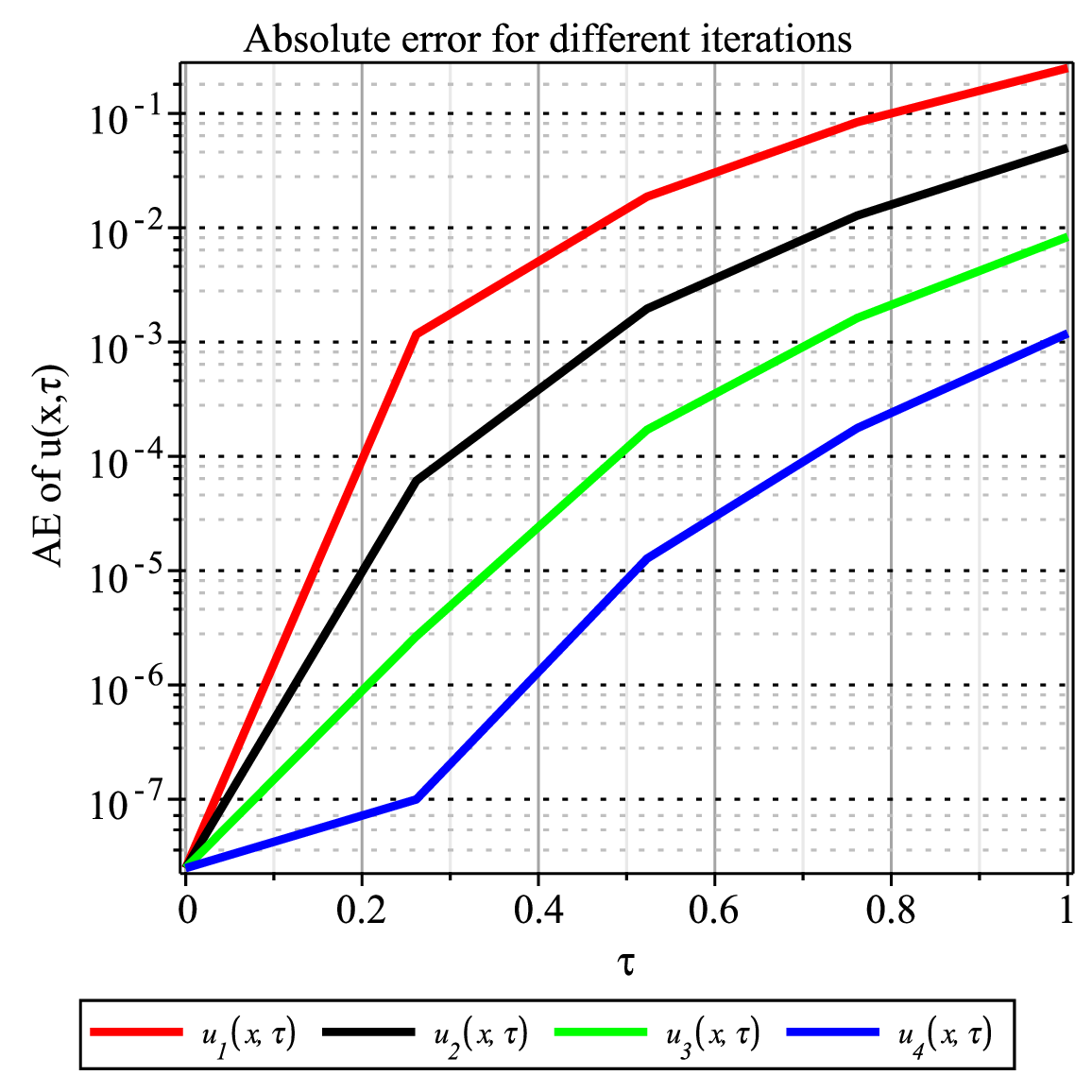}\hspace{0.1em}
    \caption{  {Comparison plots of Exact, Classical, and Modified LADM solutions, along with their corresponding absolute errors for Problem \ref{problem1}.}}\label{fig1}
  \end{figure}
\begin{figure}[H]
    \centering
     \includegraphics[width=5.45cm]{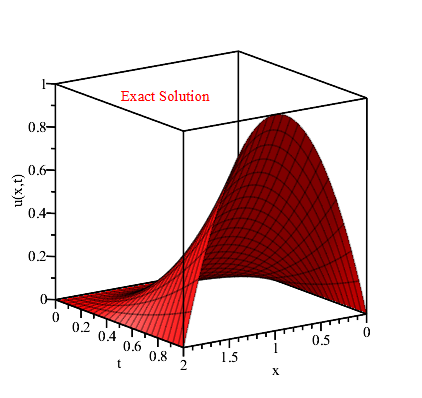}\hspace{0.00em}%
     \includegraphics[width=5.45cm]{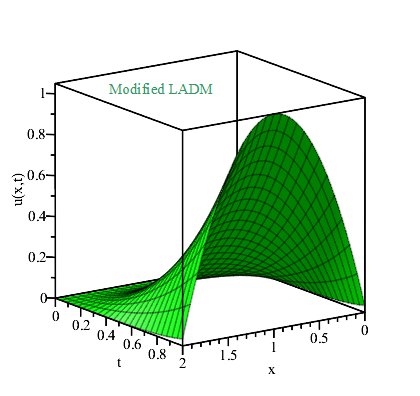}\hspace{0.00em}
     \includegraphics[width=5.45cm]{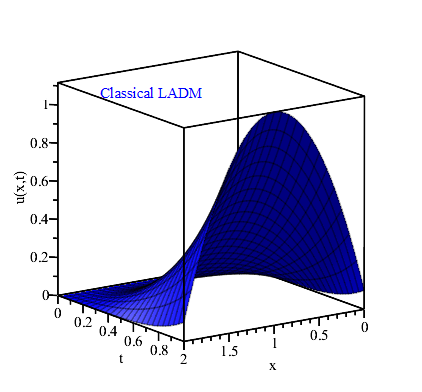}\hspace{0.000em}
    \caption{3D comparison plot presenting the Exact, Classical, and Modified LADM solutions for  problem \ref{problem1}.}\label{fig2}
  \end{figure}
\end{prb}
\begin{prb}\label{problem2} Consider the following two dimensional fractional order PDEs
\begin{equation}\label{EQQ1}
\frac{\partial^{\alpha} u(x, y, t)}{\partial t^{\alpha}} + u(x, y, t) = \frac{\partial^{2} u}{\partial x^{2}} + \frac{\partial^{2} u}{\partial y^{2}} + h(x, y, t), \quad 0 \leq x, y \leq 2, \quad 0 < \alpha \leq 1, \quad t \geq 0,
\end{equation}
with the IBCs as follows:
\begin{equation}\nonumber
\begin{split}
u(x, y, 0) &= 0, \\
u(0, y, t) &= u(2, y, t) = t^{2} y (2 - y), \\
u(x, 0, t) &= u(x, 2, t) = t^{2} x (2 - x),
\end{split}
\end{equation}
where the source term $h(x,y,t)$ is given by
\begin{equation}\nonumber
h(x, y, t) = \frac{2}{\Gamma(3 - \alpha)} t^{2 - \alpha} \left[ x(2 - x) + y(2 - y) \right] + t^{2} \left[ x(2 - x) + y(2 - y) \right] + 4t^{2}.
\end{equation}
The exact solution  for  \ref{problem2} is given by
\begin{equation}\nonumber
u(x, y, t) = t^{2} \left[ x(2 - x) + y(2 - y) \right].
\end{equation}
\begin{figure}[H]
    \centering
    \includegraphics[width=5.4cm]{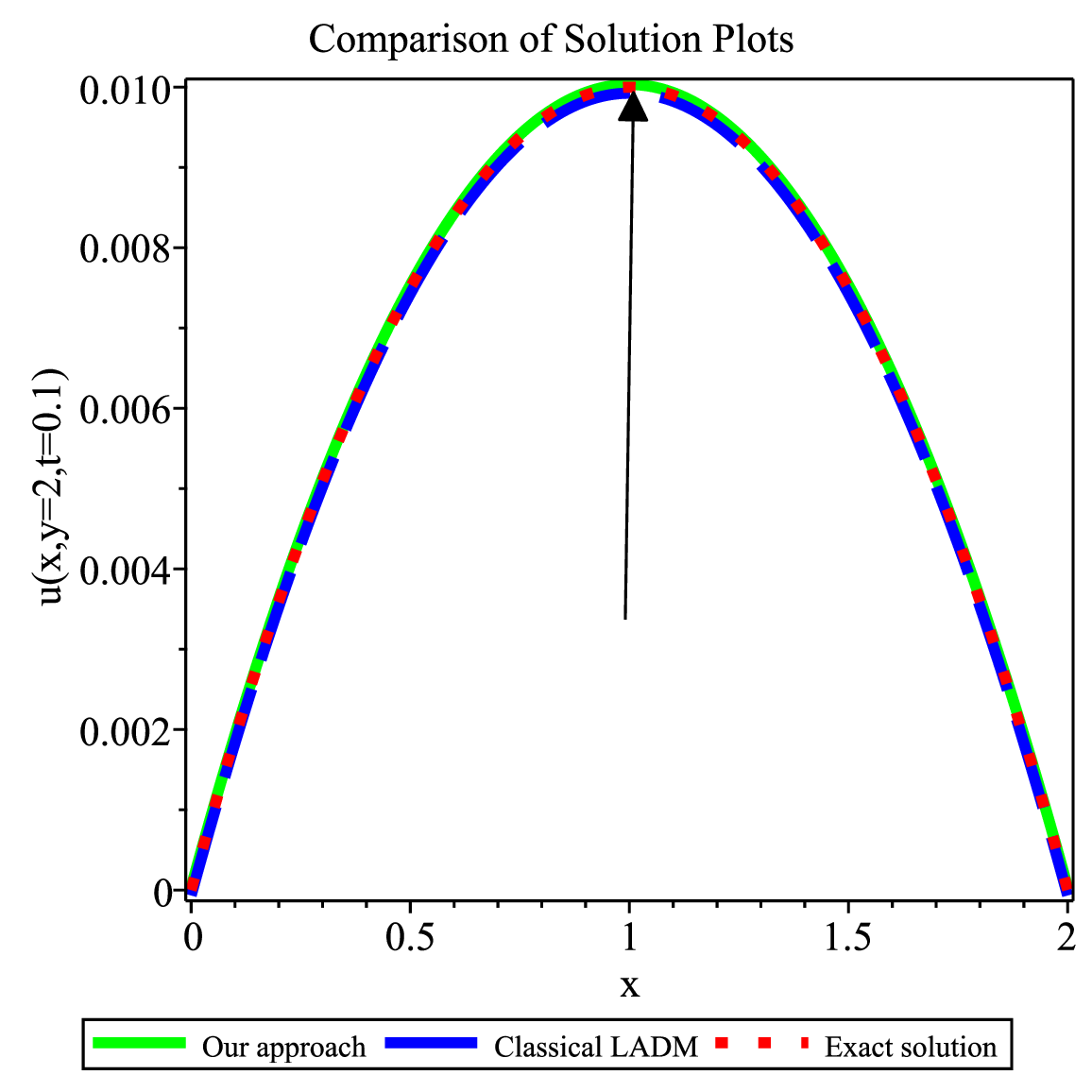}%
    \llap{\raisebox{0.95cm}{%
      \includegraphics[width=2.50cm]{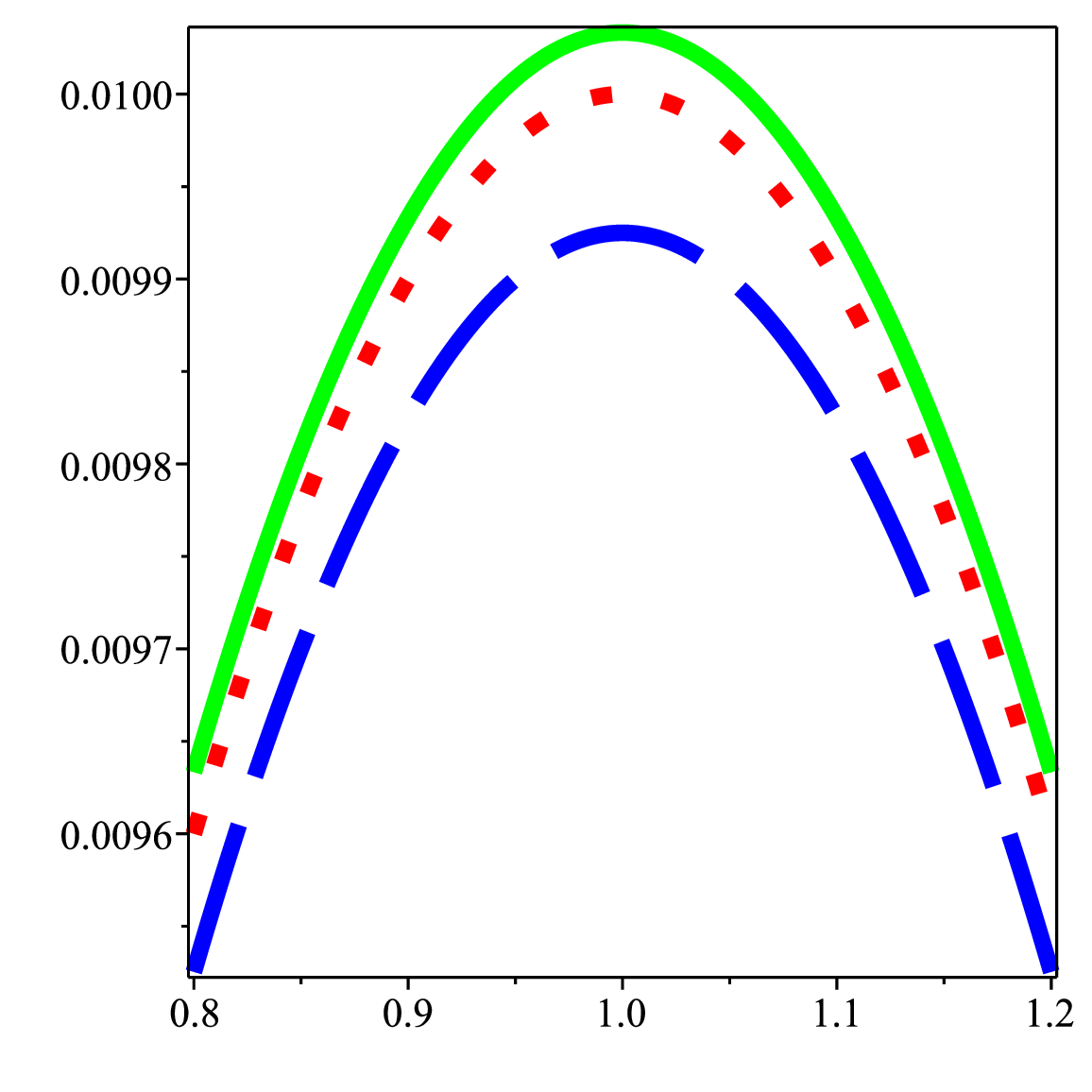}\hspace{2.6 em}%
    }}\hspace{0.1em}
     \includegraphics[width=5.4cm]{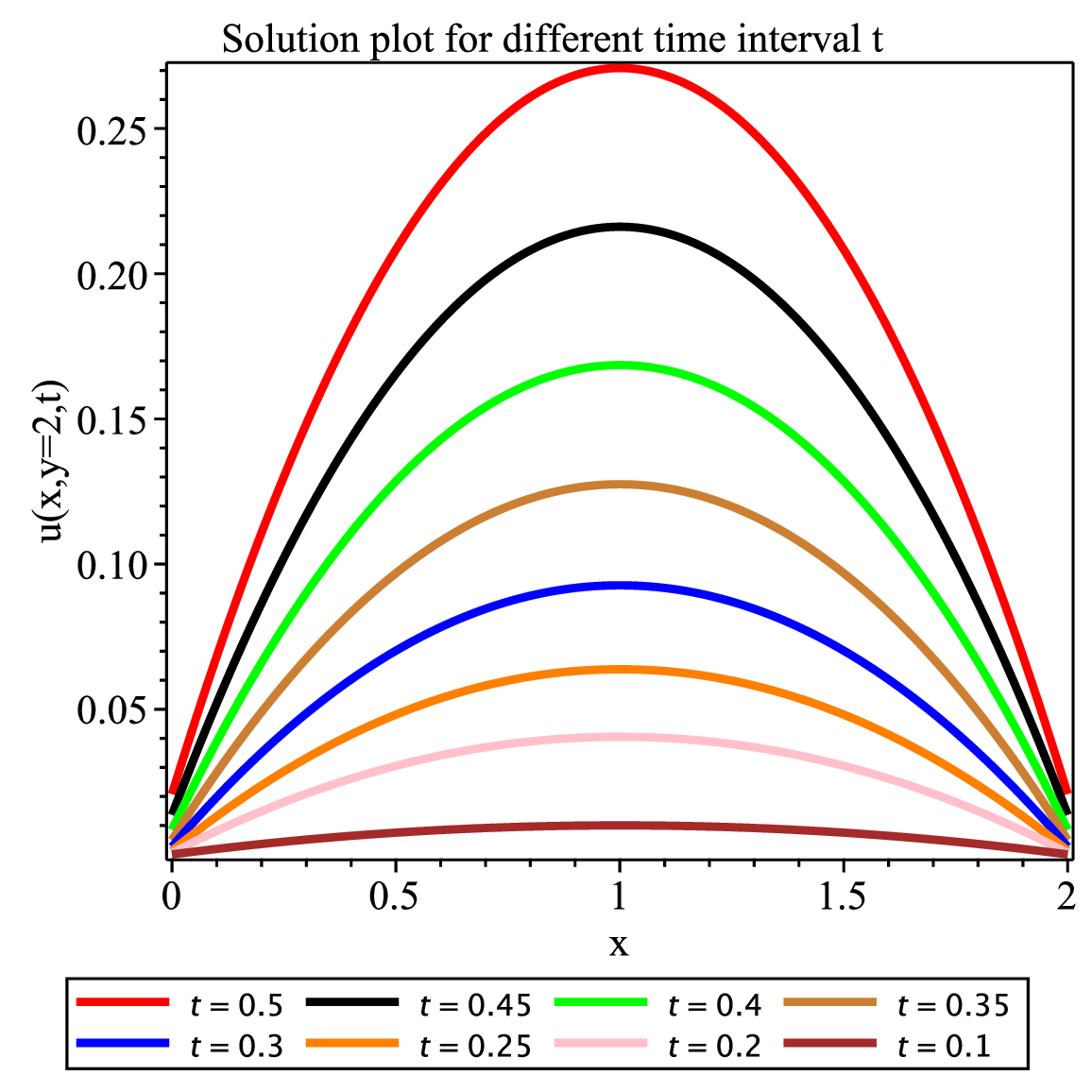}\hspace{0.1em}%
     \includegraphics[width=5.4cm]{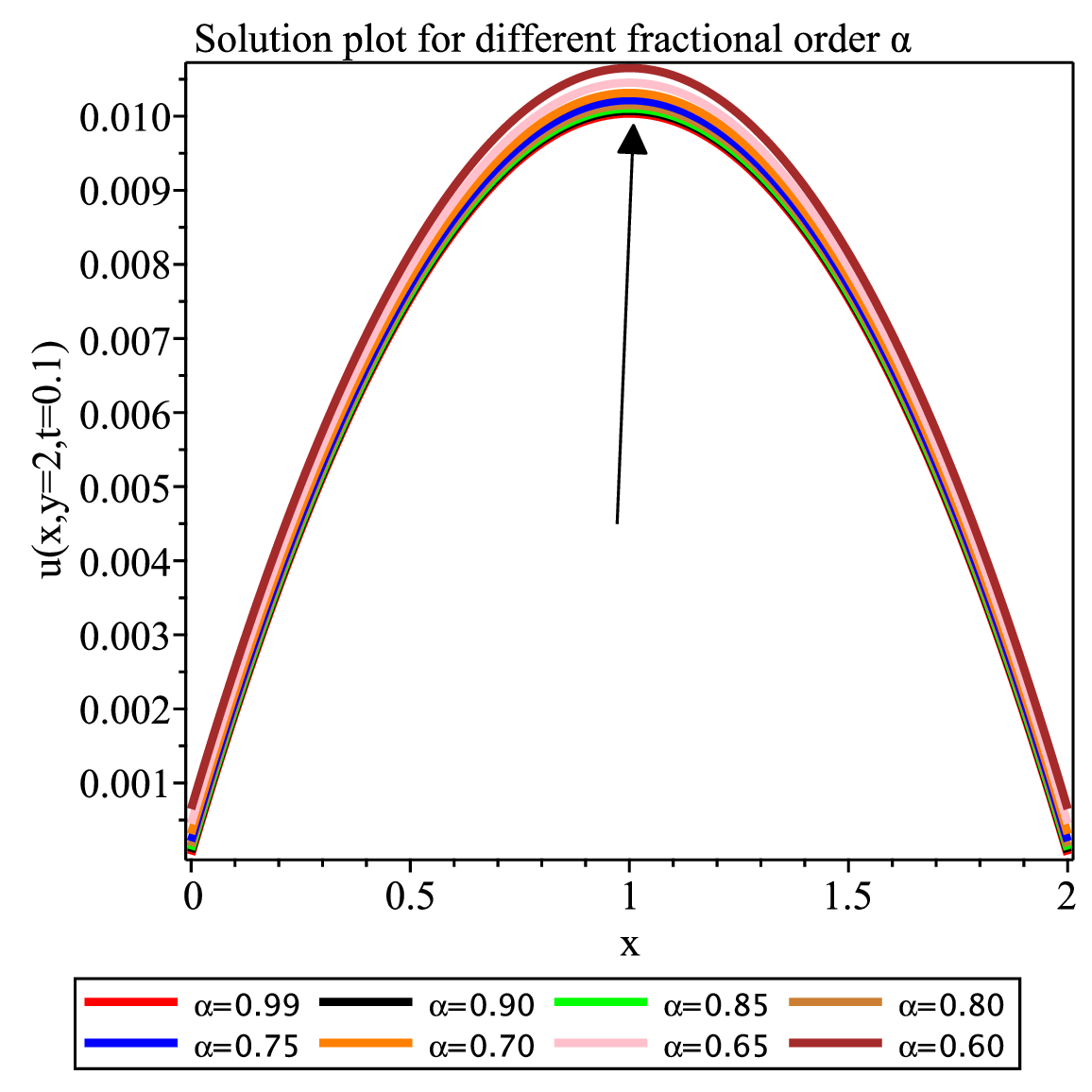}%
    \llap{\raisebox{1.25cm}{%
      \includegraphics[width=2.30cm]{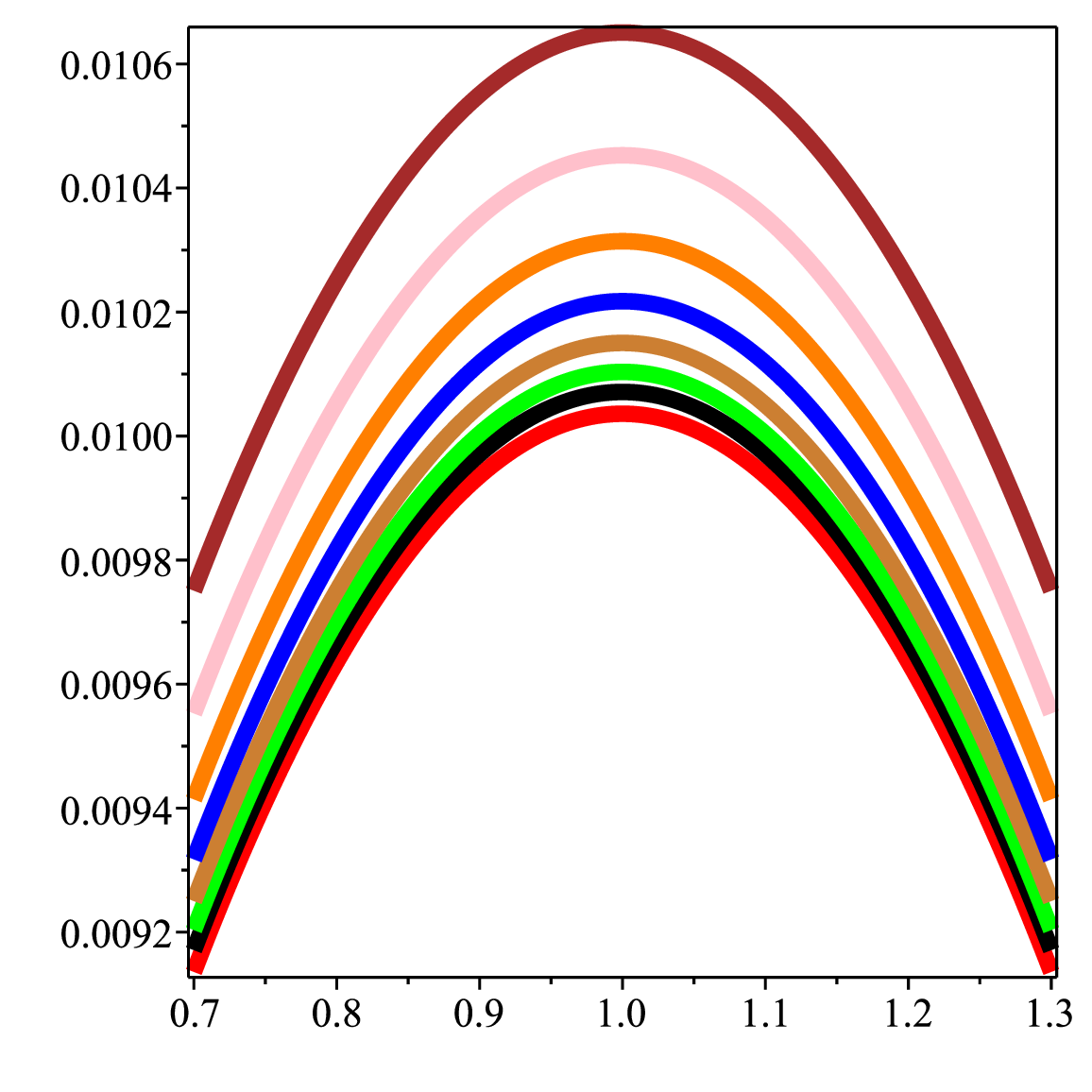}\hspace{2.90em}%
    }}\hspace{0.1em}
    \caption{  {Comparison plots of Exact, Classical, and Modified LADM solutions   at different fractional order $\alpha$ for Problem \ref{problem2}.}}\label{fig3}
  \end{figure}
\begin{figure}[H]
    \centering
     \includegraphics[width=5.45cm]{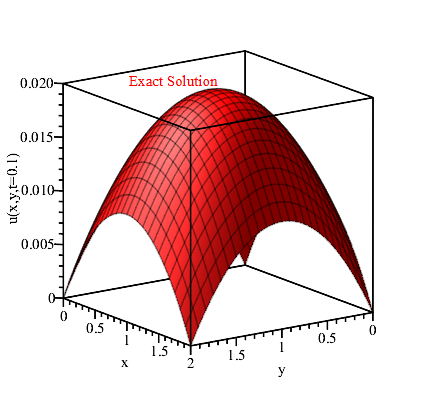}\hspace{0.00em}%
     \includegraphics[width=5.45cm]{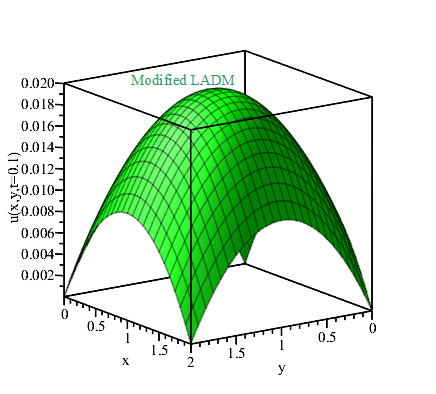}\hspace{0.00em}
     \includegraphics[width=5.45cm]{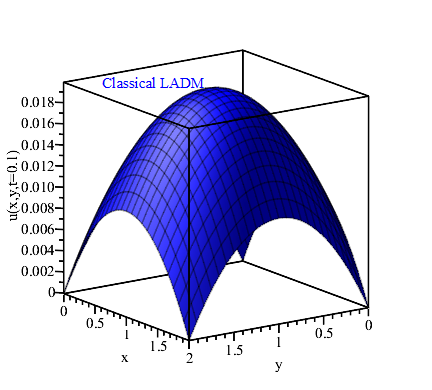}\hspace{0.000em}
    \caption{3D comparison plot presenting the Exact, Classical, and Modified LADM solutions for  problem \ref{problem2}.}\label{fig4}
  \end{figure}
\end{prb}
\begin{prb}\label{problem32}
Consider the following initial and boundary value problem
\begin{equation}
\frac{\partial^\alpha u(x,t)}{\partial t^\alpha}  =  \frac{\partial u(x,t)}{\partial x}+h({{{x}}}, {{{t}}}), \quad t > 0, \quad 0 \leq x \leq 1, \quad 0 < \alpha \leq 1,
\end{equation}
with the IBCs condition as follows: 
\begin{equation*}
\begin{split} 
u(x, 0) &= {{\rm e}^{x}},\\
u(0, t) &= {t}^{3}+1, \\
u(1, t) &=    {t}^{3}\cos(1)+ {\rm e^{1}},
\end{split}
\end{equation*}
where $h(x,t)$ is given by
\begin{equation*}
h({{{x}}}, {{{t}}})=\left( {\frac {6{t}^{3-\alpha}}{\Gamma \left( 4-\alpha \right) }}+{
t}^{3} \right)\cos(x)-  {{\rm e}^{x}}
.
\end{equation*}
The exact solution  for  \ref{problem32} is given by
\begin{equation*}
u(x, t) =   {t}^{3}\cos(x)+ { \rm e}^{x}.
\end{equation*}

\begin{figure}[H]
    \centering
    \includegraphics[width=5.4cm]{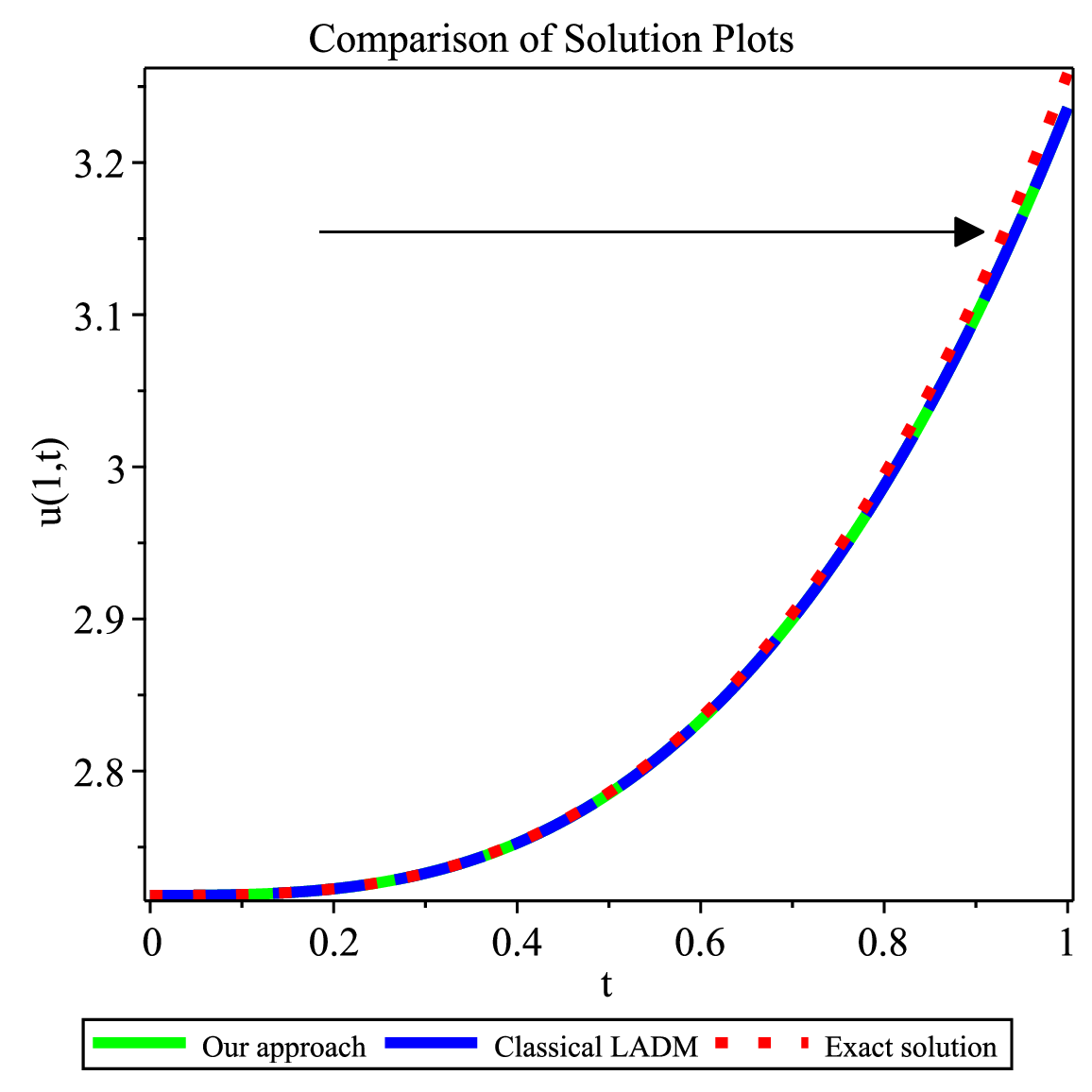}%
    \llap{\raisebox{2.15cm}{%
      \includegraphics[width=2.90cm]{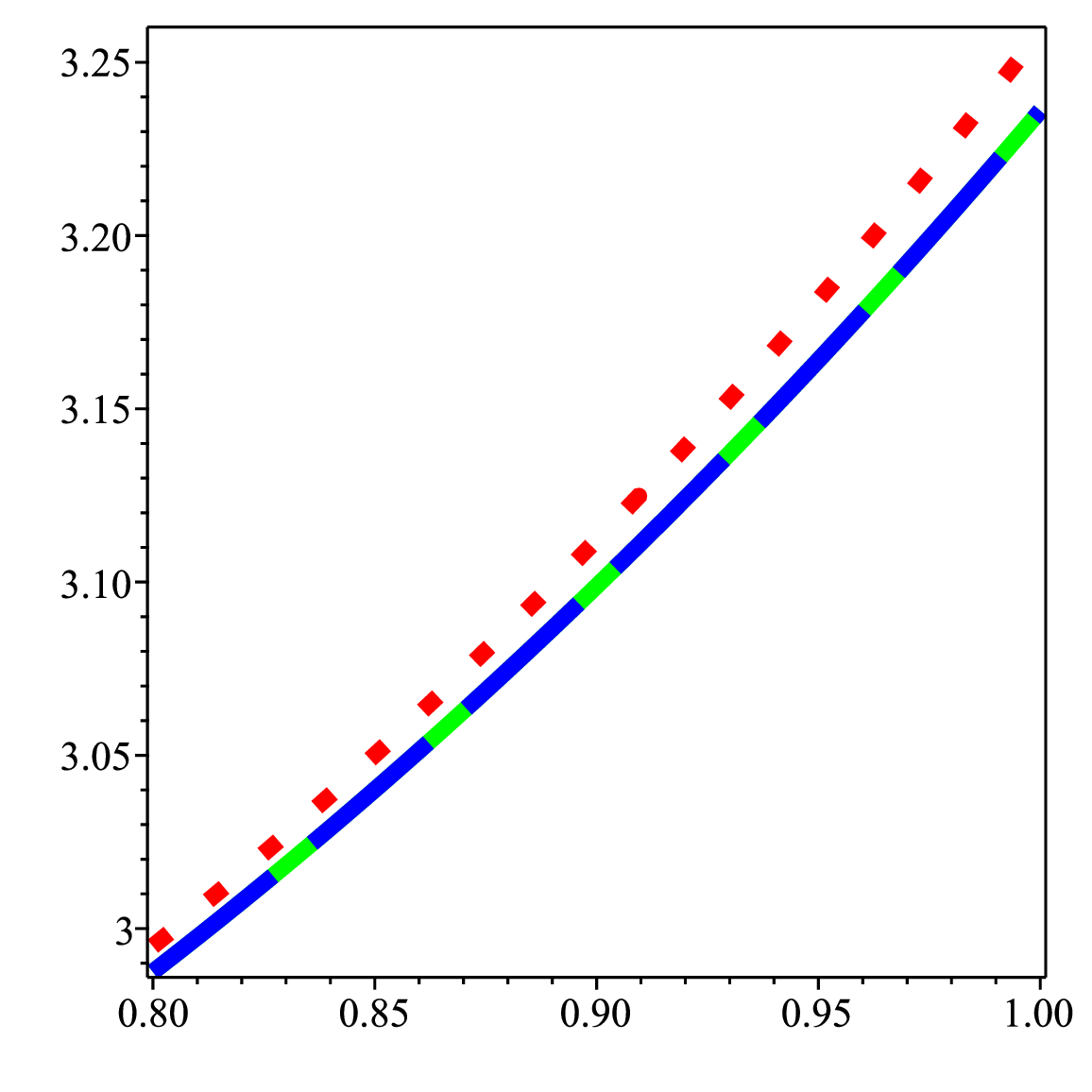}\hspace{4.25em}%
    }}\hspace{0.1em}
     \includegraphics[width=5.4cm]{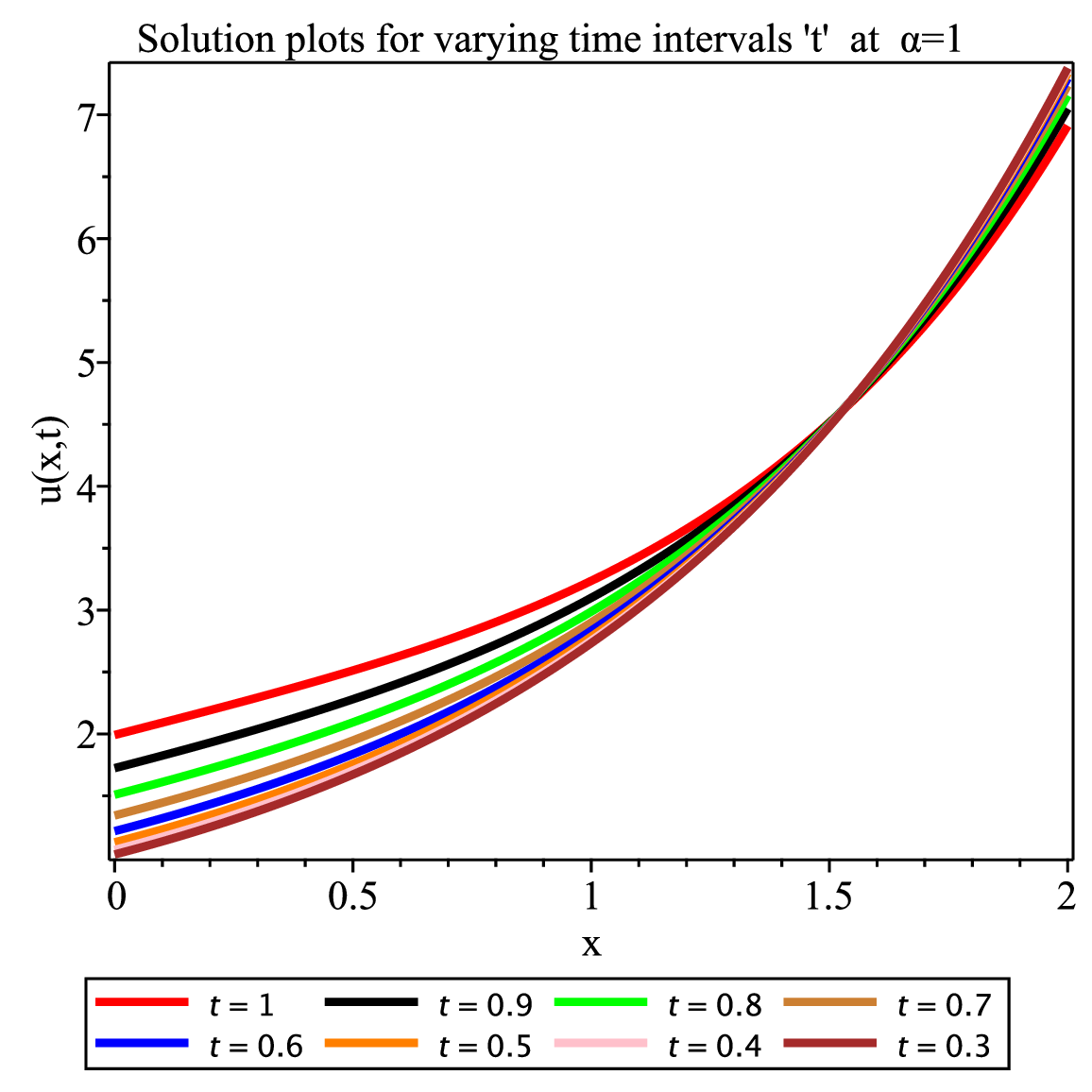}\hspace{0.1em}%
     \includegraphics[width=5.4cm]{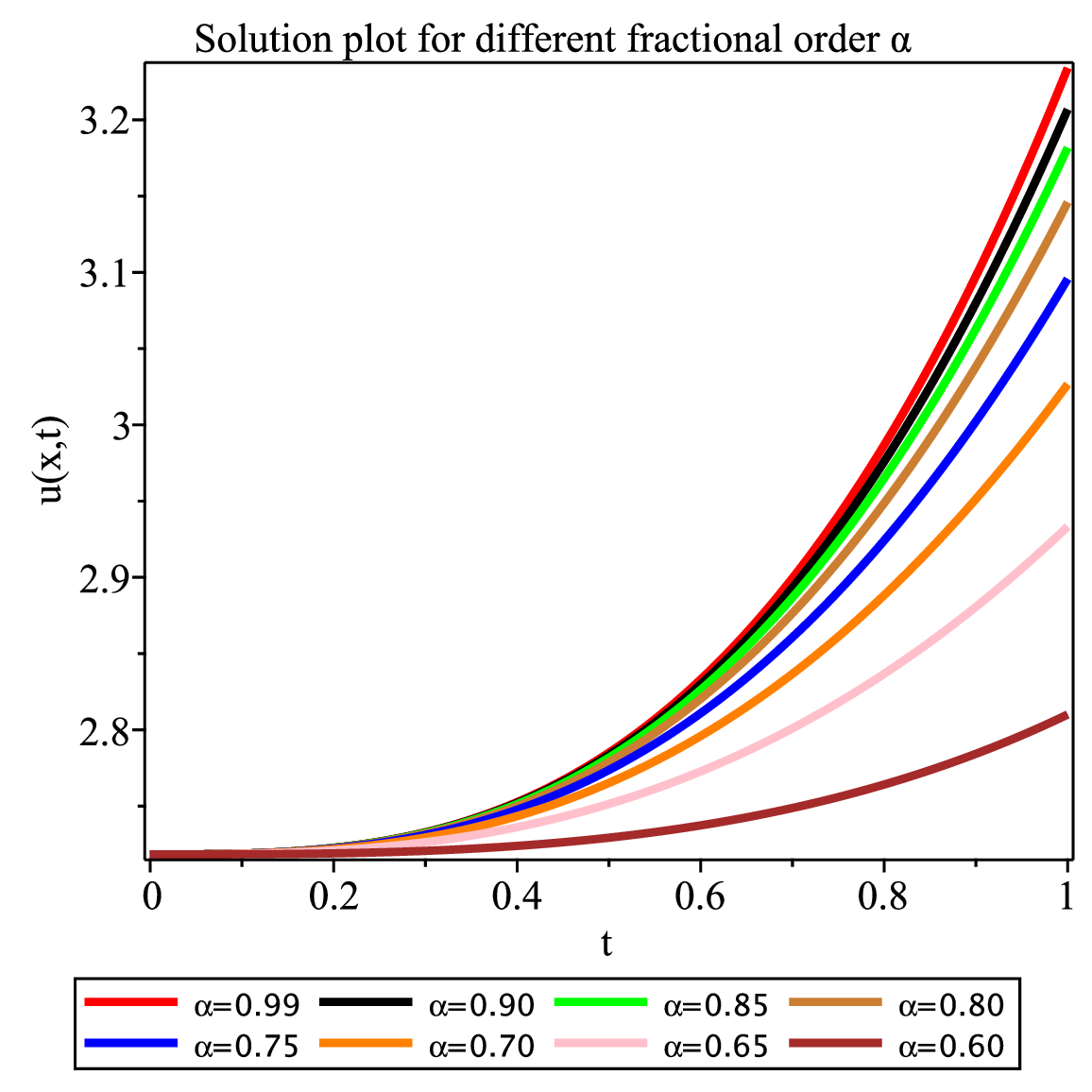}\hspace{0.1em}
    \caption{  {Comparison plots of Exact, Classical, and Modified LADM solutions   at different fractional order $\alpha$ for Problem \ref{problem32}.}}\label{fig5}
  \end{figure}
\begin{figure}[H]
    \centering
     \includegraphics[width=5.45cm]{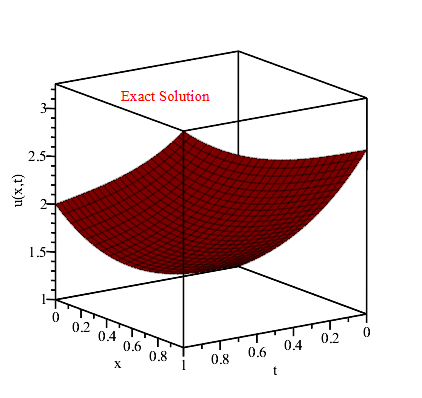}\hspace{0.00em}%
     \includegraphics[width=5.45cm]{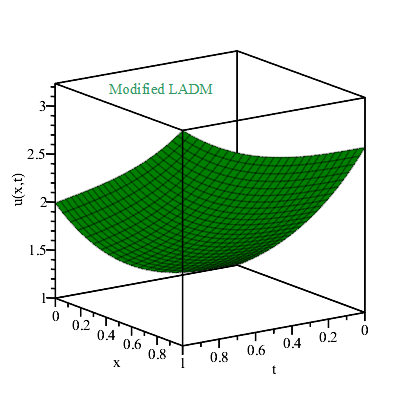}\hspace{0.00em}
     \includegraphics[width=5.45cm]{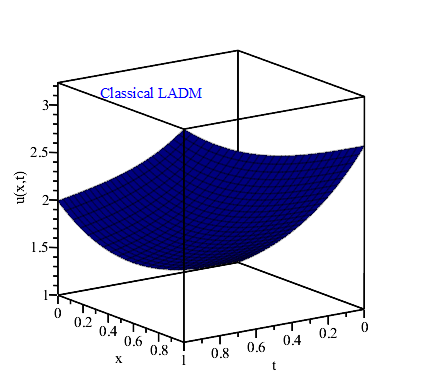}\hspace{0.000em}
    \caption{3D comparison plot presenting the Exact, Classical, and Modified LADM solutions for  problem \ref{problem32}.}\label{fig6}
  \end{figure}

\end{prb}
\begin{prb}\label{newp3}
Consider the following initial and boundary value problem:
\begin{equation}
\frac{\partial^\alpha u(x,t)}{\partial t^\alpha} +  \frac{\partial u(x,t)}{\partial x} + \frac{\partial^2 u(x,t)}{\partial x^2} = h({{{x}}}, {{{t}}}), \quad t \geq 0, \quad 0 \leq x \leq 1, \quad 0 < \alpha \leq 1,
\end{equation}
with the IBCs condition as follows:
\begin{equation*}
\begin{split}
u(x, 0) &= x^2,\\
u(0, t) &=  1, \\
u(1, t) &=    {t}^{3+{\it \alpha}}\sin \left( 1 \right) +1,
\end{split}
\end{equation*}
where $h(x,t)$ is given by
\begin{equation*}
h({{{x}}}, {{{t}}})= \left( 1/6\,\Gamma \left( 4+\alpha \right) {t}^{3}+{t}^{3+\alpha}
 \right) \sin \left( x \right)
.
\end{equation*}
The exact solution  for  \ref{newp3} is given by
\begin{equation*}
u(x, t) ={t}^{3+\alpha}\sin \left( x \right) +1
.
\end{equation*}
\begin{figure}[H]
    \centering
    \includegraphics[width=5.4cm]{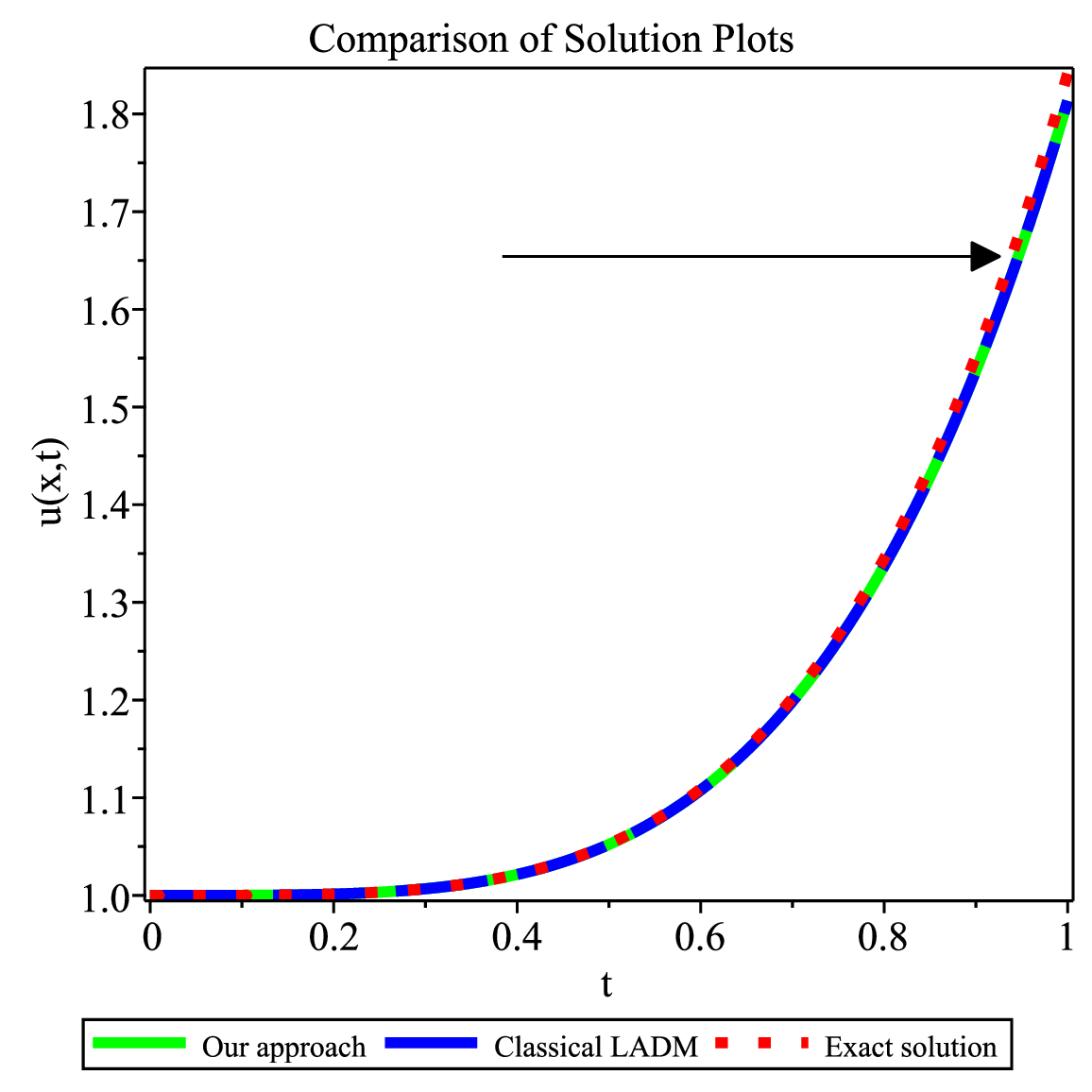}%
    \llap{\raisebox{1.895cm}{%
      \includegraphics[width=3.150cm]{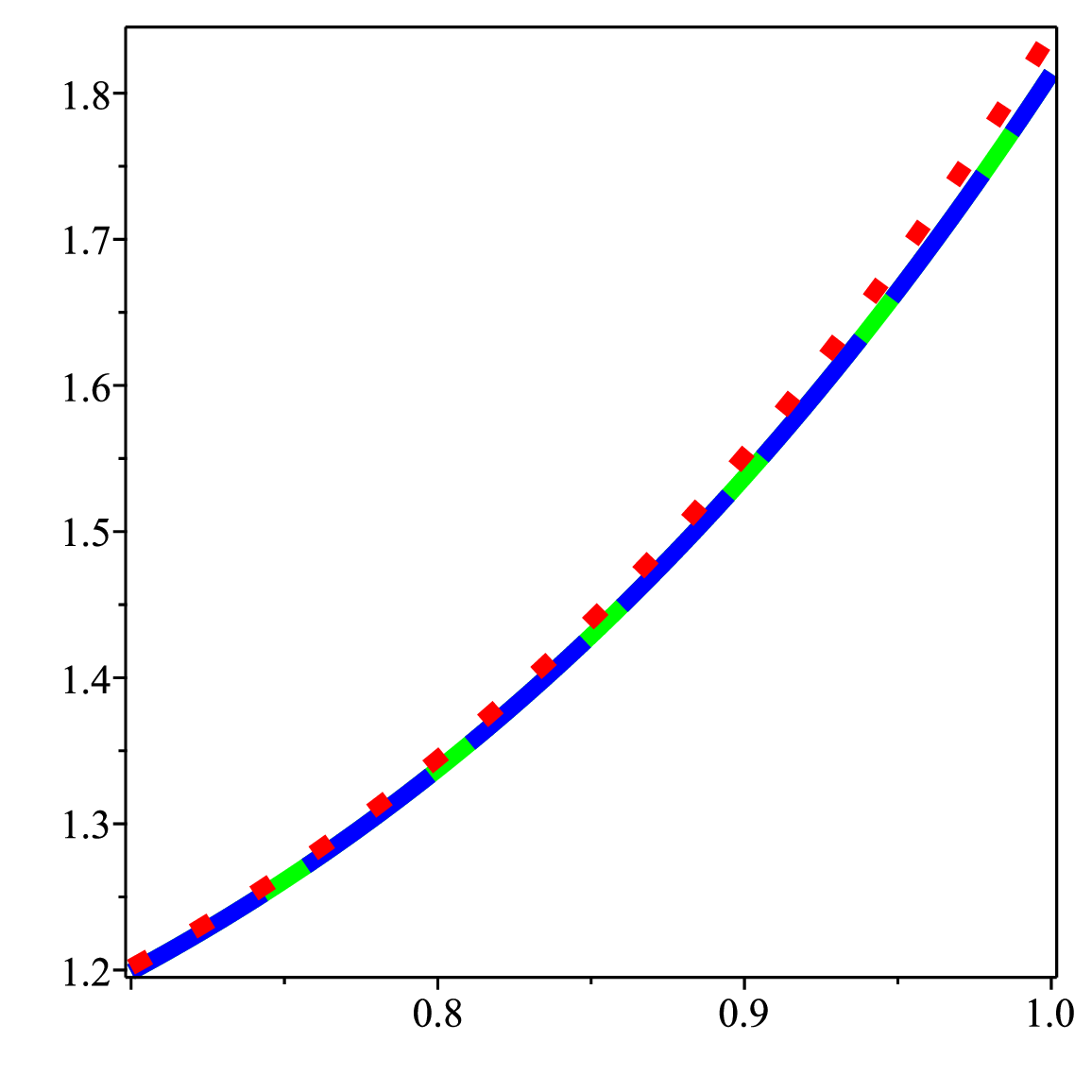}\hspace{3.69715em}%
    }}\hspace{0.1em}
     \includegraphics[width=5.4cm]{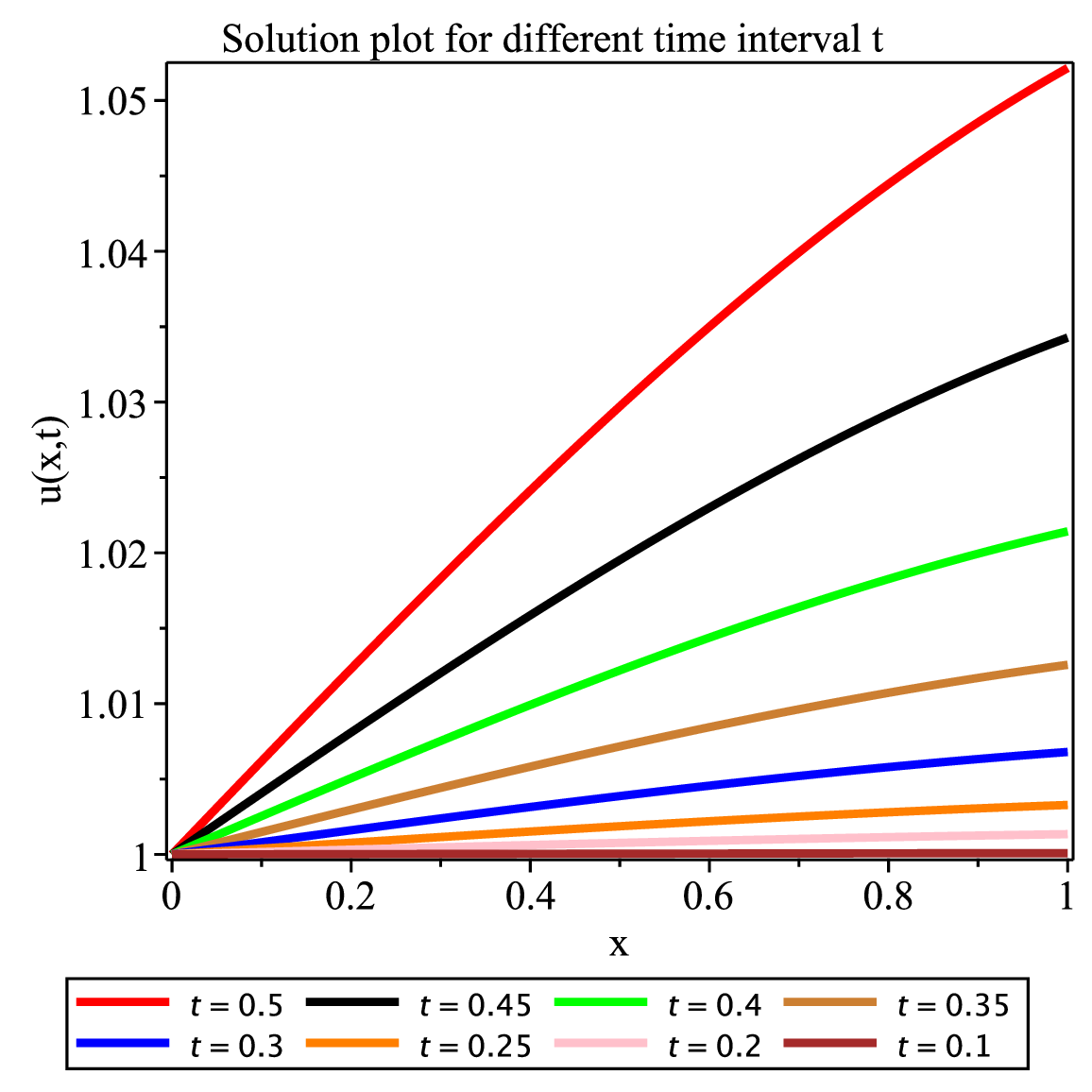}\hspace{0.1em}%
     \includegraphics[width=5.4cm]{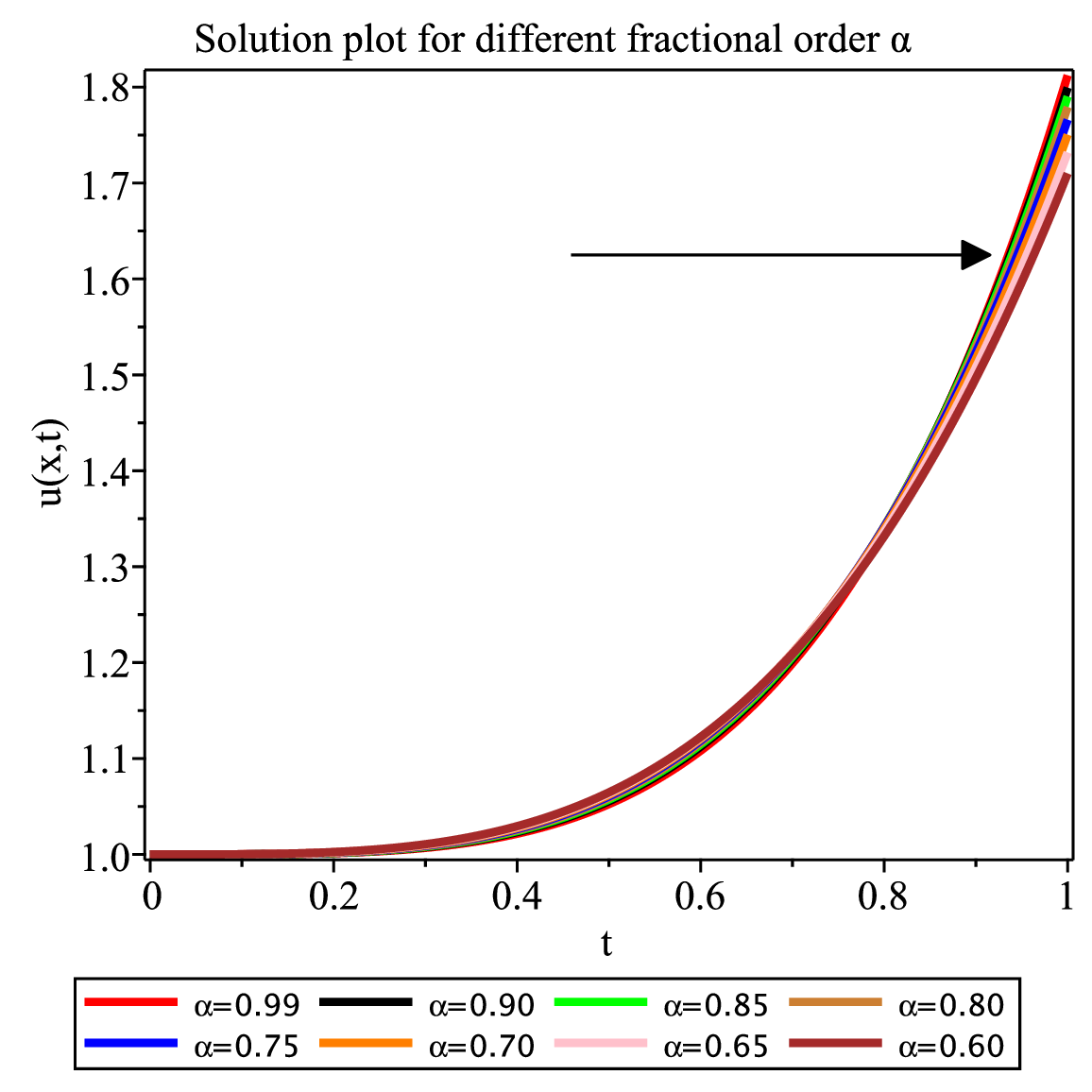}%
    \llap{\raisebox{1.895cm}{%
      \includegraphics[width=2.850cm]{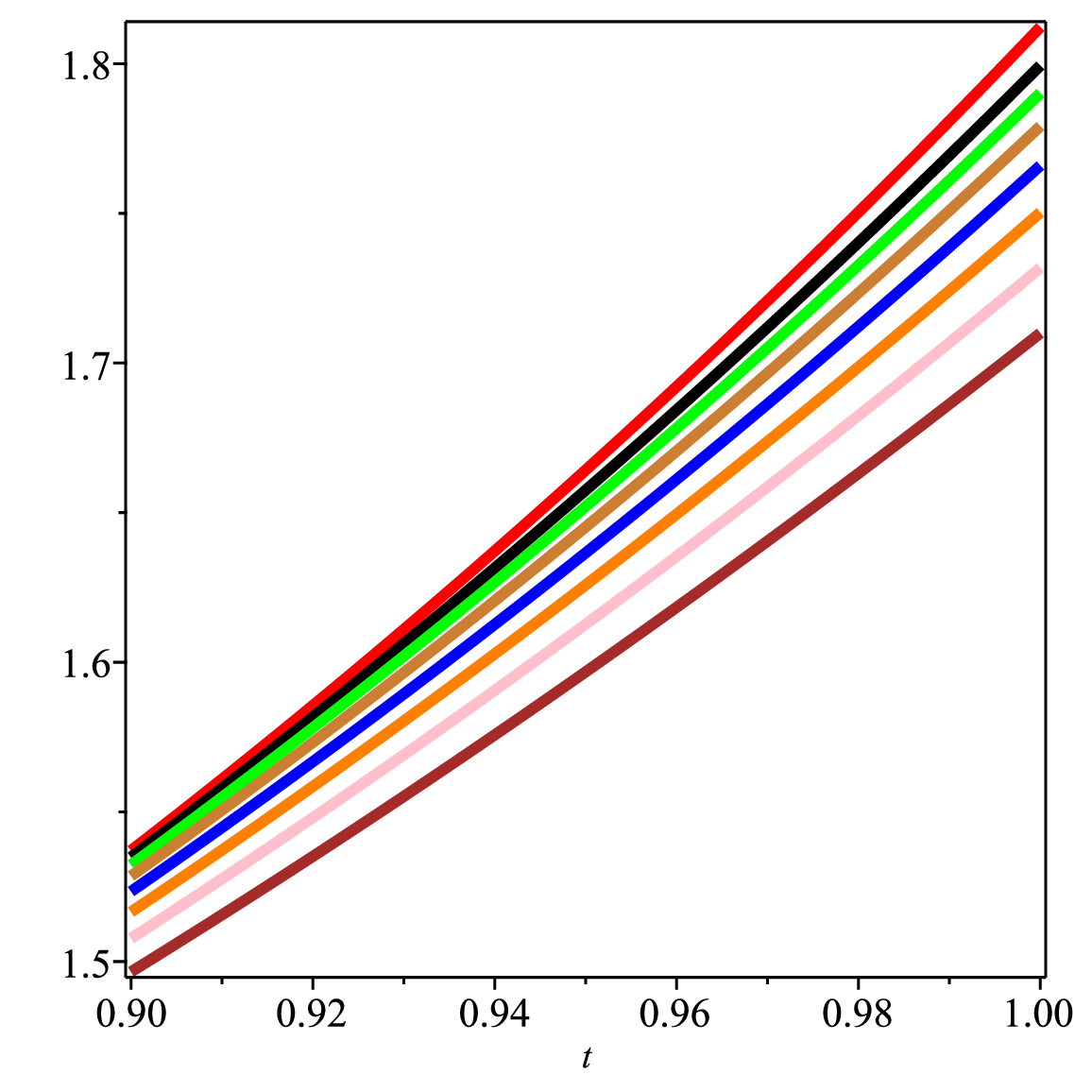}\hspace{4.339715em}%
    }}\hspace{0.1em}
    \caption{  {Comparison plots of Exact, Classical, and Modified LADM solutions   at different fractional order $\alpha$ for Problem \ref{newp3}.}}\label{fig7}
  \end{figure}

\begin{figure}[H]
    \centering
     \includegraphics[width=5.45cm]{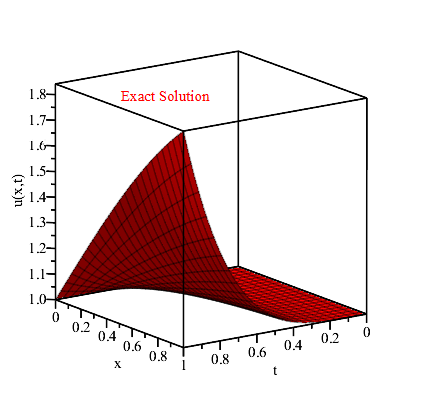}\hspace{0.00em}%
     \includegraphics[width=5.45cm]{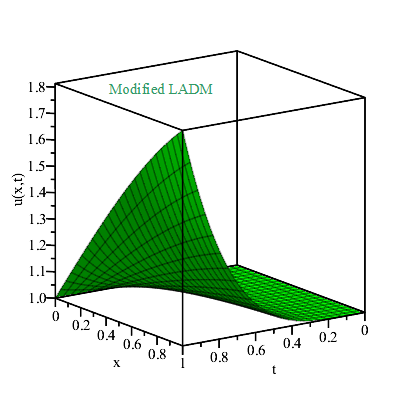}\hspace{0.00em}
     \includegraphics[width=5.45cm]{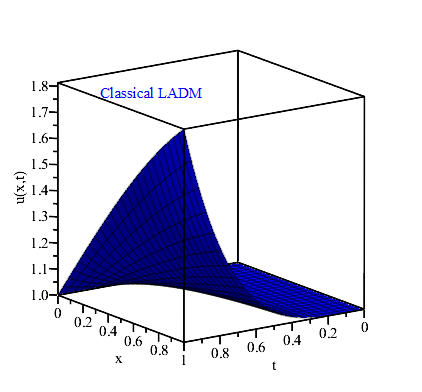}\hspace{0.000em}
    \caption{3D comparison plot presenting the Exact, Classical, and Modified LADM solutions for  problem \ref{newp3}.}\label{fig8}
  \end{figure}
\end{prb}
\begin{prb}\label{problem3}
Consider the following initial and boundary value problem
\begin{equation}
\frac{\partial^\alpha u(x,t)}{\partial t^\alpha} + \frac{\partial u(x,t)}{\partial x} = h({{{x}}}, {{{t}}}), \quad t > 0, \quad 0 \leq x \leq 1, \quad 0 < \alpha \leq 1,
\end{equation}
with the IBCs condition as follows:
\begin{equation*}
\begin{split}
u(x, 0) &= 0,\\
u(0, t) &=  0, \\
u(1, t) &=    t \sin(1),
\end{split}
\end{equation*}
where $h(x,t)$ is given by
\begin{equation*}
h({{{x}}}, {{{t}}})=\frac{{t}^{1-\alpha} \sin(x)}{\Gamma(2 - \alpha)} + t \cos(x).
\end{equation*}
The exact solution  for  \ref{problem3} is given by
\begin{equation*}
u(x, t) = t \sin(x).
\end{equation*}

\begin{figure}[H]
    \centering
    \includegraphics[width=5.4cm]{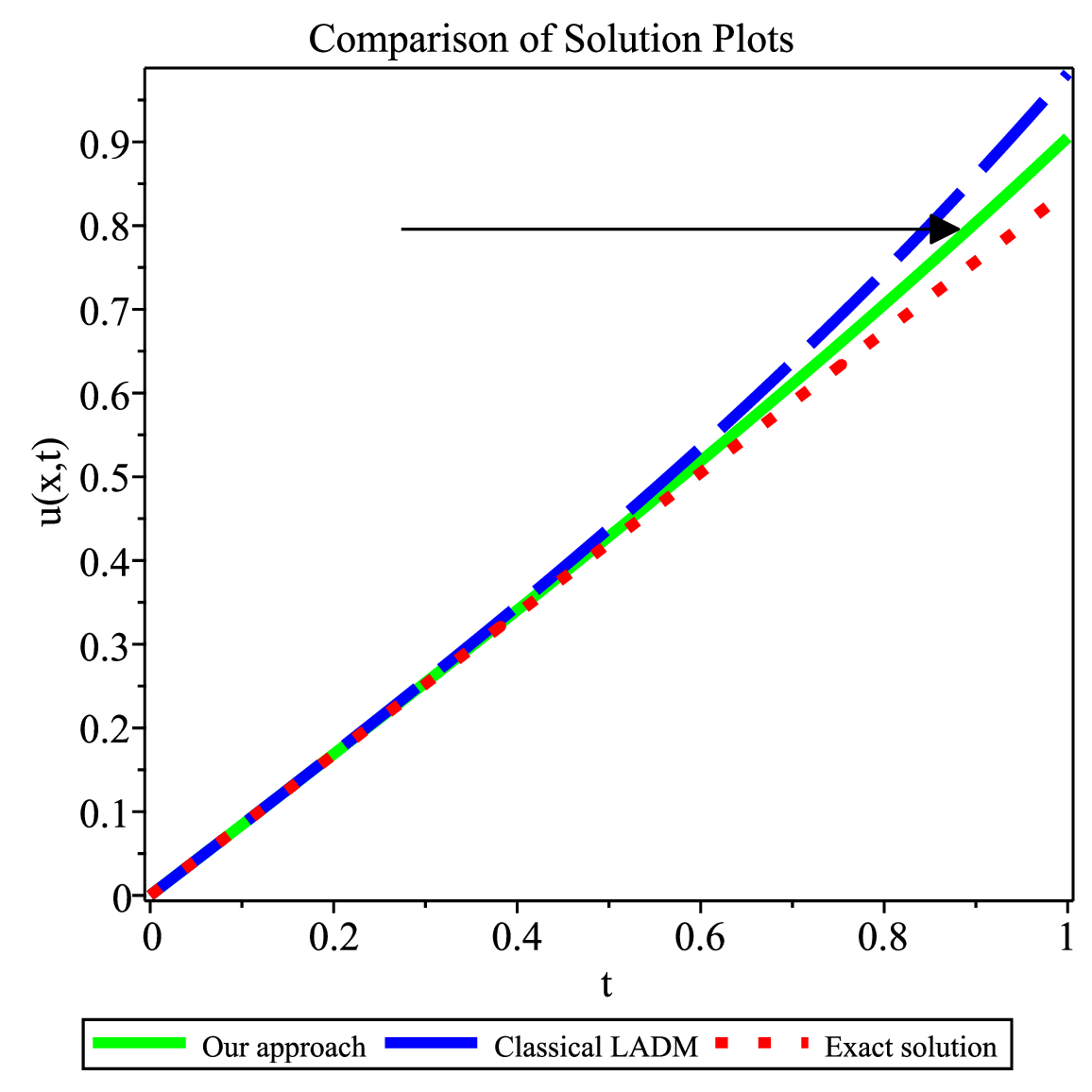}%
    \llap{\raisebox{2.85cm}{%
      \includegraphics[width=2.20cm]{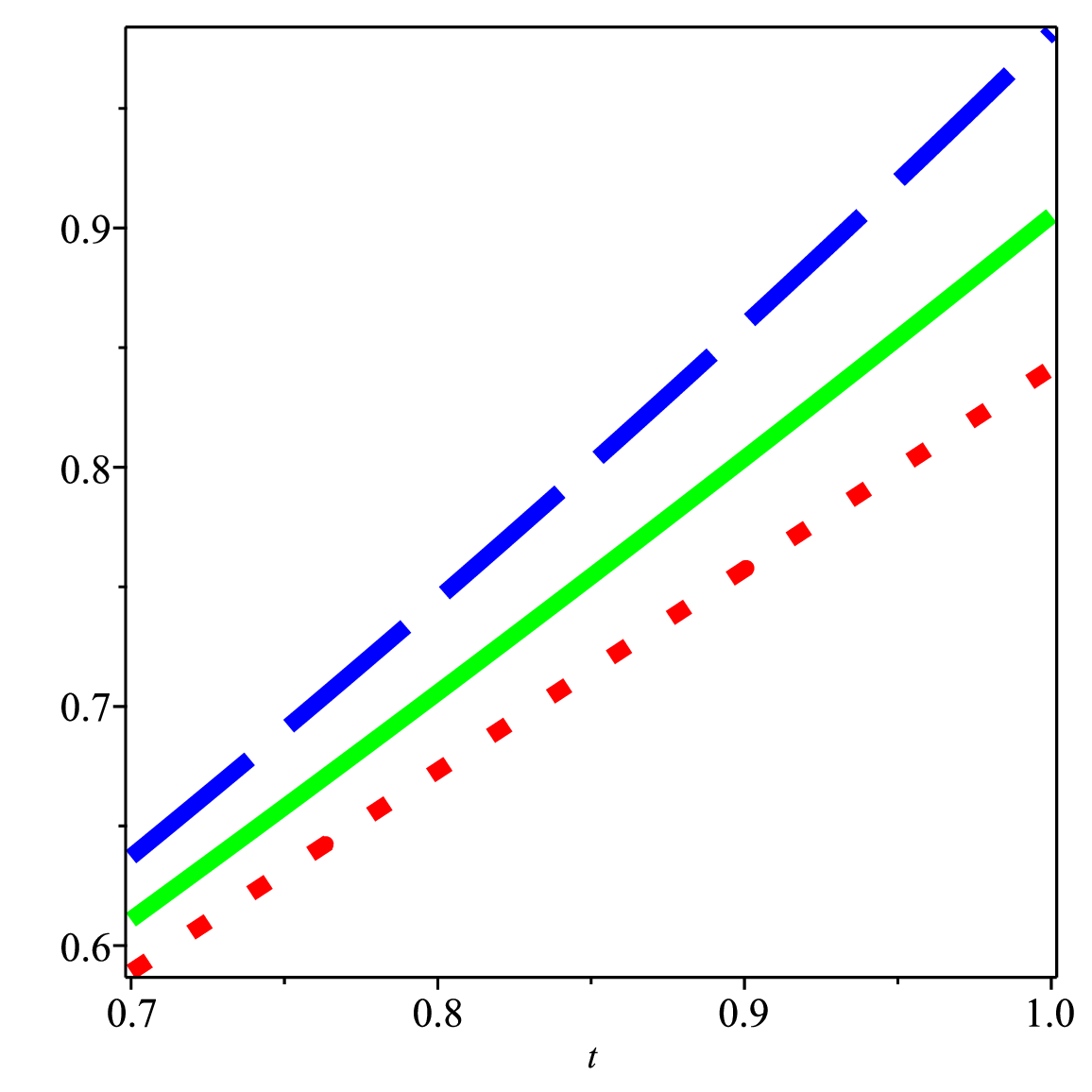}\hspace{5.85em}%
    }}\hspace{0.1em}
     \includegraphics[width=5.4cm]{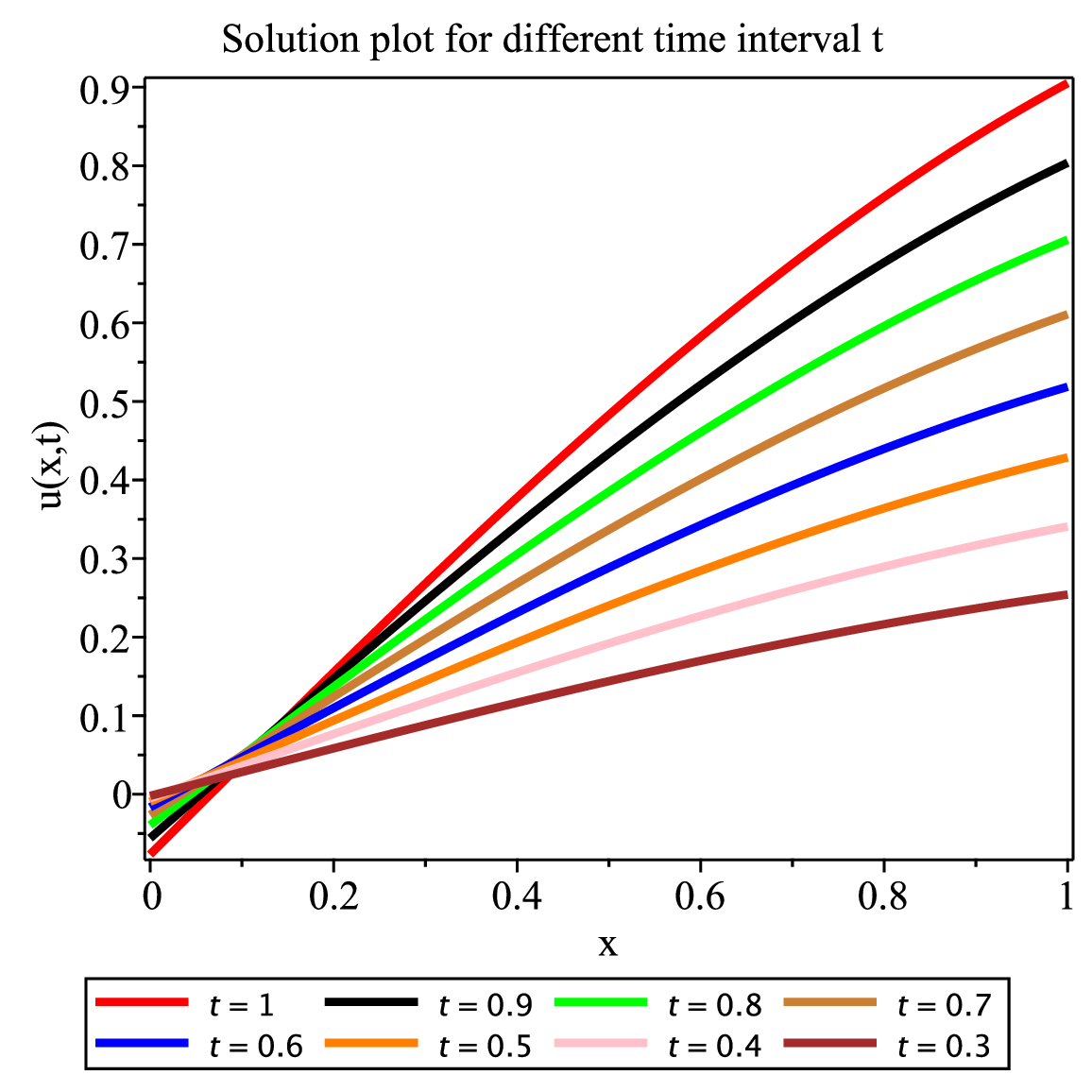}\hspace{0.1em}%
     \includegraphics[width=5.4cm]{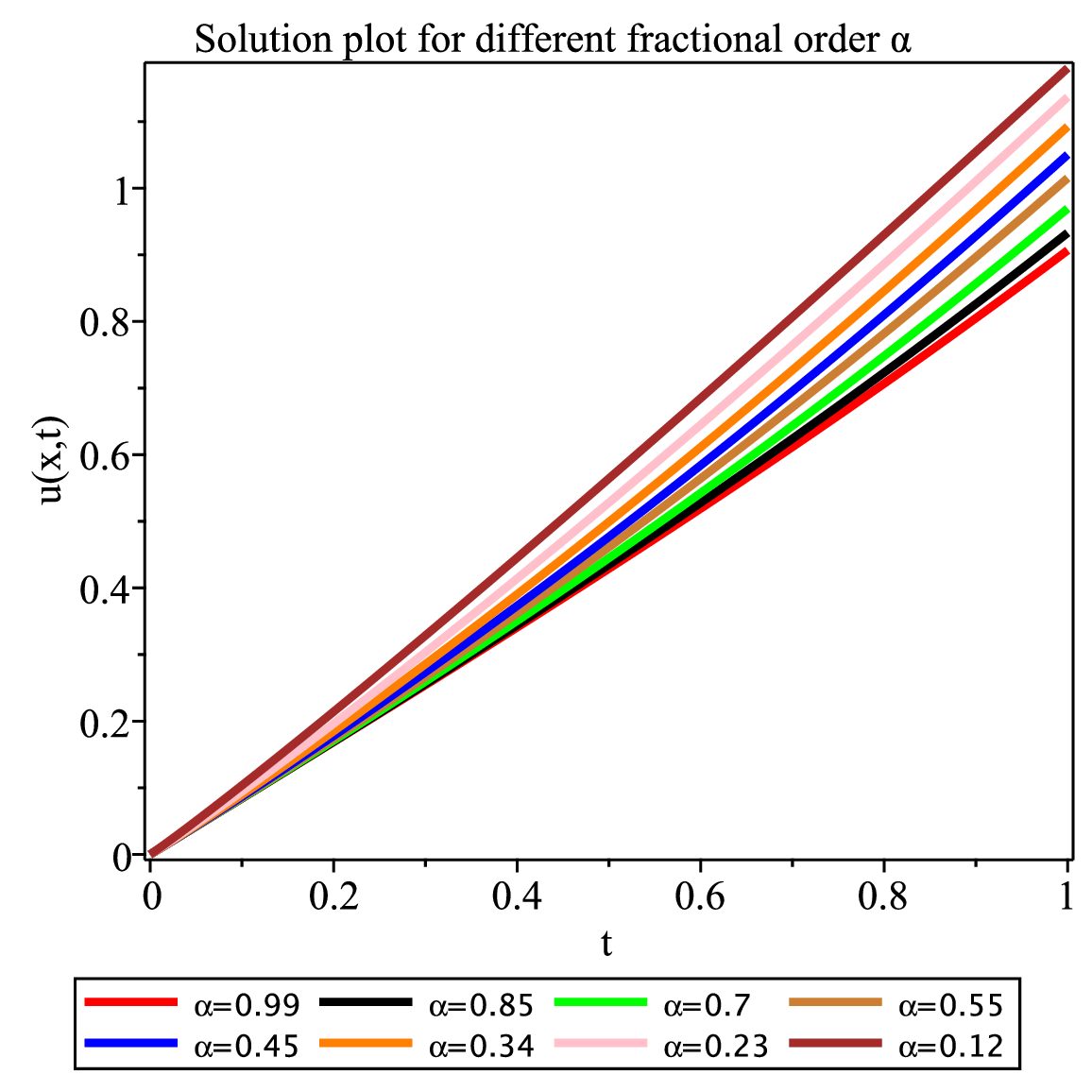}\hspace{0.1em}
    \caption{  {Comparison plots of Exact, Classical, and Modified LADM solutions   at different fractional order $\alpha$ for Problem \ref{problem3}.}}\label{fig9}
  \end{figure}
\begin{figure}[H]
    \centering
     \includegraphics[width=5.45cm]{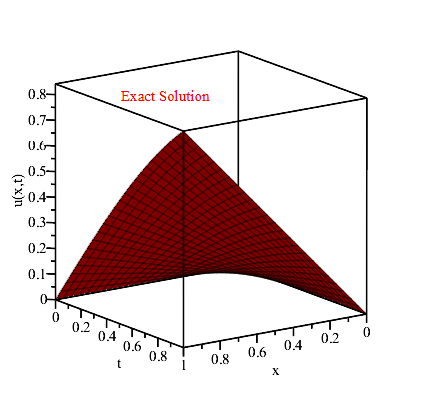}\hspace{0.00em}%
     \includegraphics[width=5.45cm]{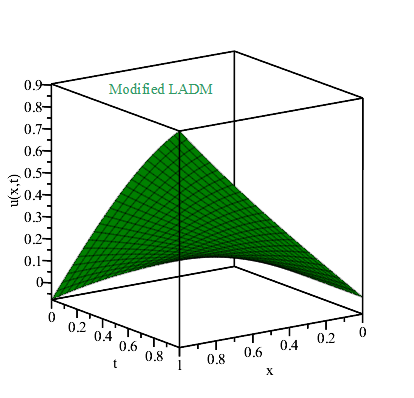}\hspace{0.00em}
     \includegraphics[width=5.45cm]{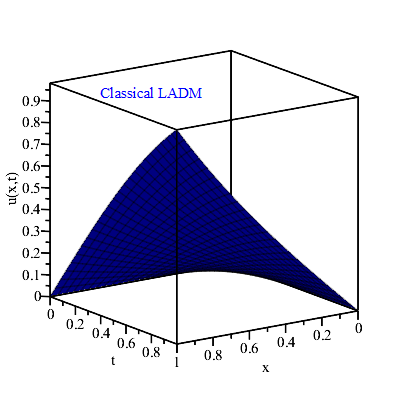}\hspace{0.000em}
    \caption{3D comparison plot presenting the Exact, Classical, and Modified LADM solutions for  problem \ref{problem3}.}\label{fig10}
  \end{figure} 
\end{prb} 
\begin{prb}\label{problem4} 
Consider the following non-linear initial and boundary value problem
\begin{equation}
\frac{\partial^\alpha u(x,t)}{\partial t^\alpha} + u^2(x, t) = h({{{x}}}, {{{t}}}), \quad t > 0, \quad 0 \leq x \leq 1, \quad 0 < \alpha \leq 1,
\end{equation}
with the IBCs condition as follows:
\begin{equation*}
\begin{split}
u(x, 0) &= 0,\\
u(0, t) &=  0, \\
u(1, t) &=    t^2,
\end{split}
\end{equation*}
 where $h(x,t)$ is given by
\begin{equation*}
h({{{x}}}, {{{t}}})=x^2( \frac{ 2t^{2 - \alpha}}{\Gamma(3 - \alpha)} + x^2 t^4).
\end{equation*}
The exact solution  for  \ref{problem4} is given by
\begin{equation*}
u(x, t) = x^2 t^2.
\end{equation*}

\begin{figure}[H]
    \centering
    \includegraphics[width=5.4cm]{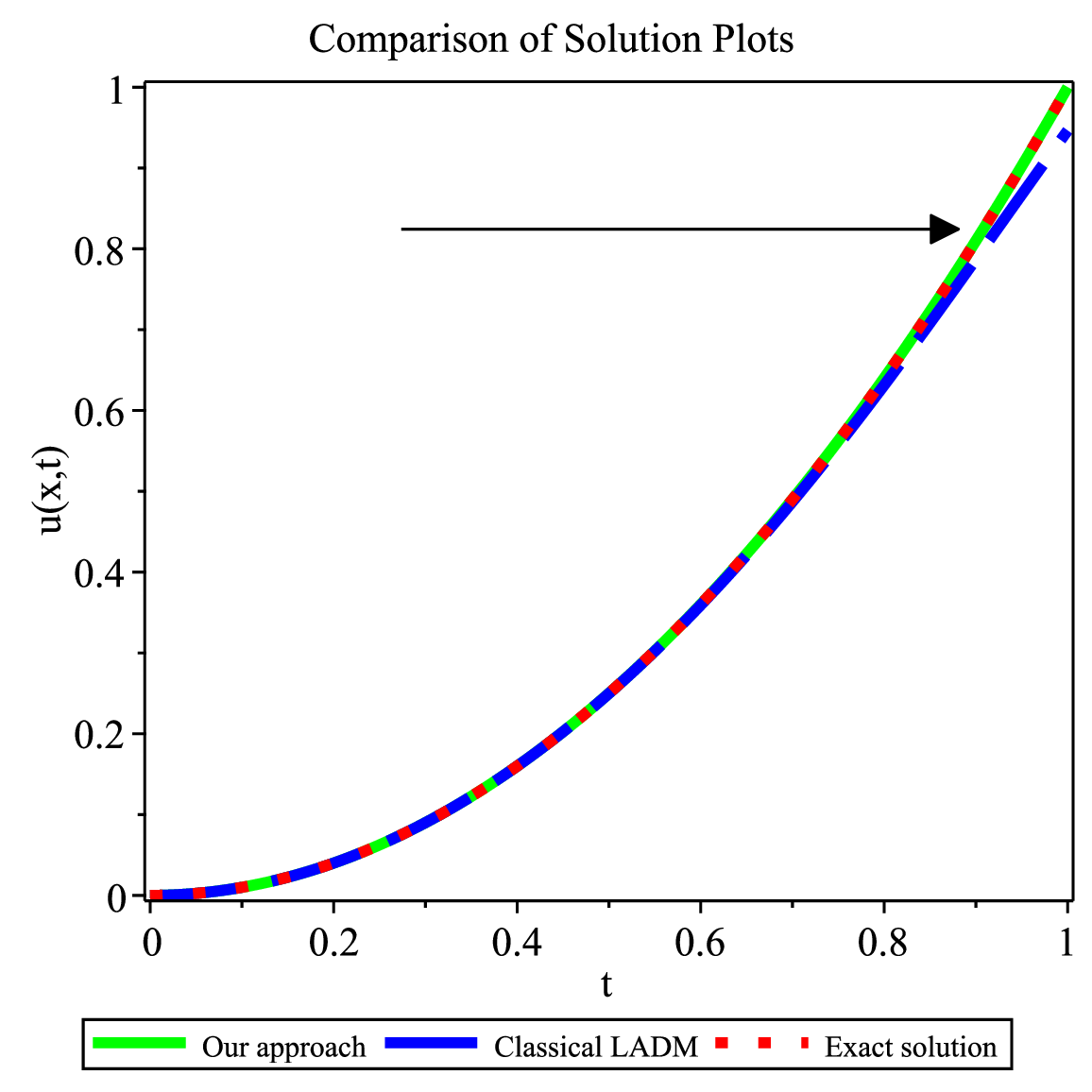}%
    \llap{\raisebox{2.45cm}{%
      \includegraphics[width=2.50cm]{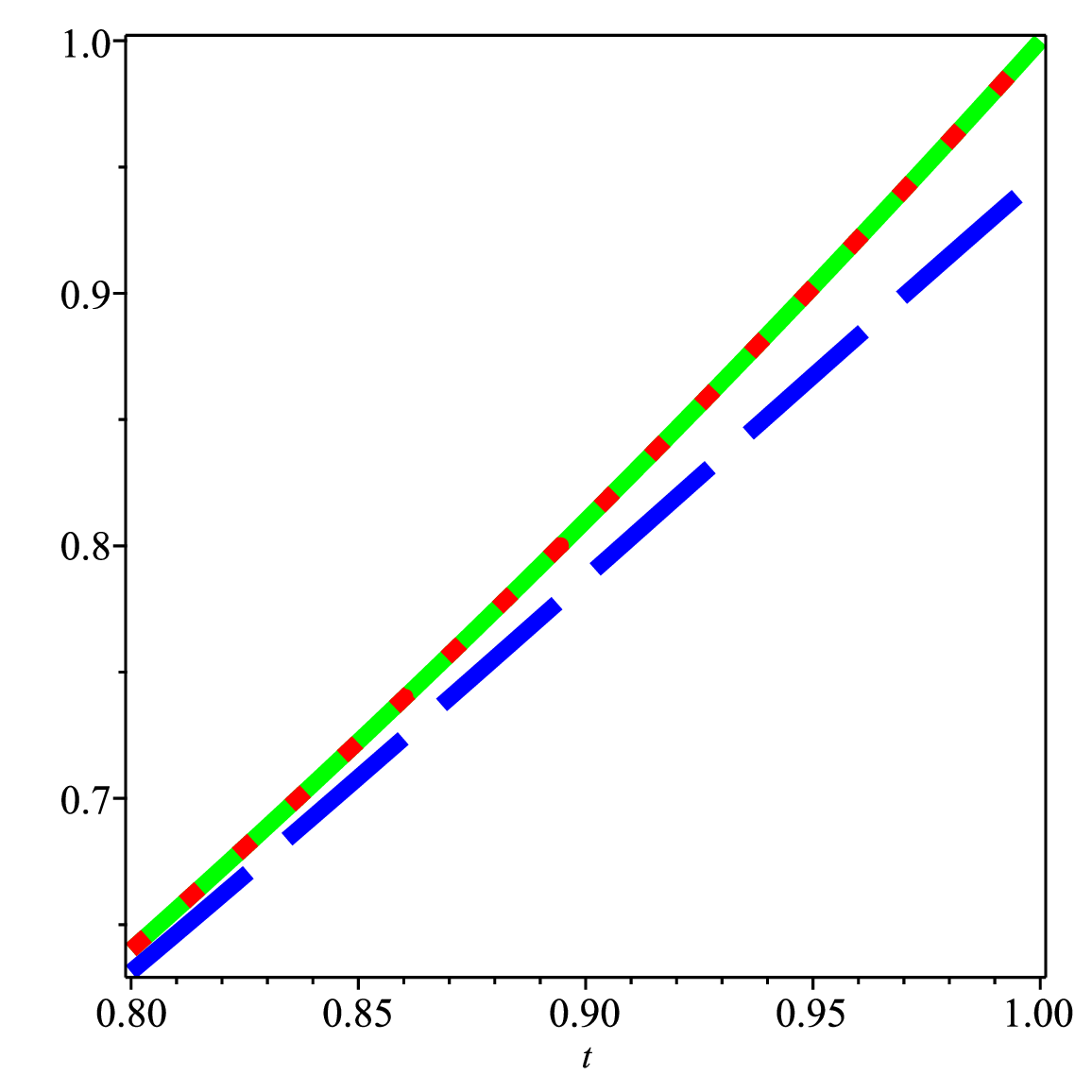}\hspace{5.05em}%
    }}\hspace{0.1em}
     \includegraphics[width=5.4cm]{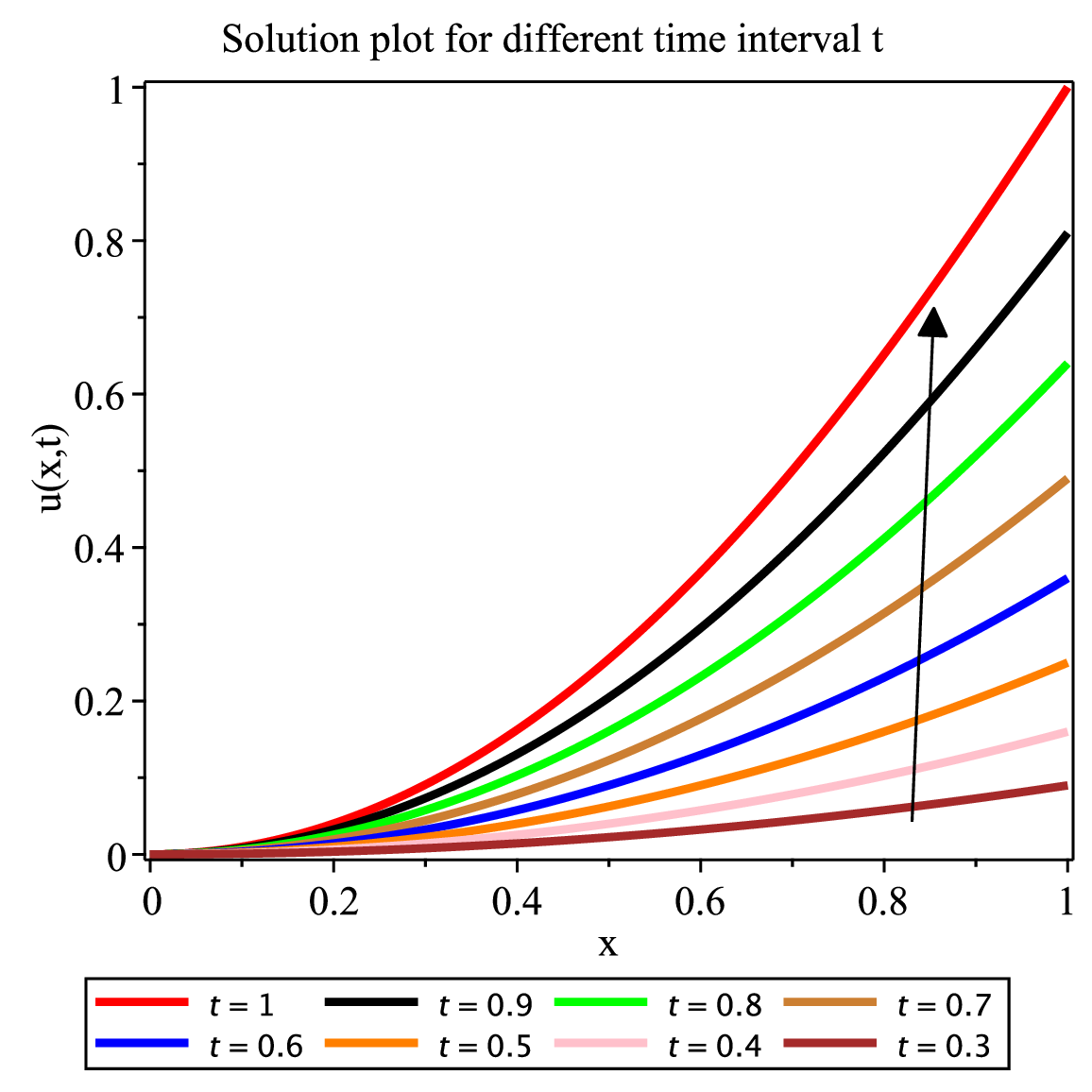}\hspace{0.1em}%
     \includegraphics[width=5.4cm]{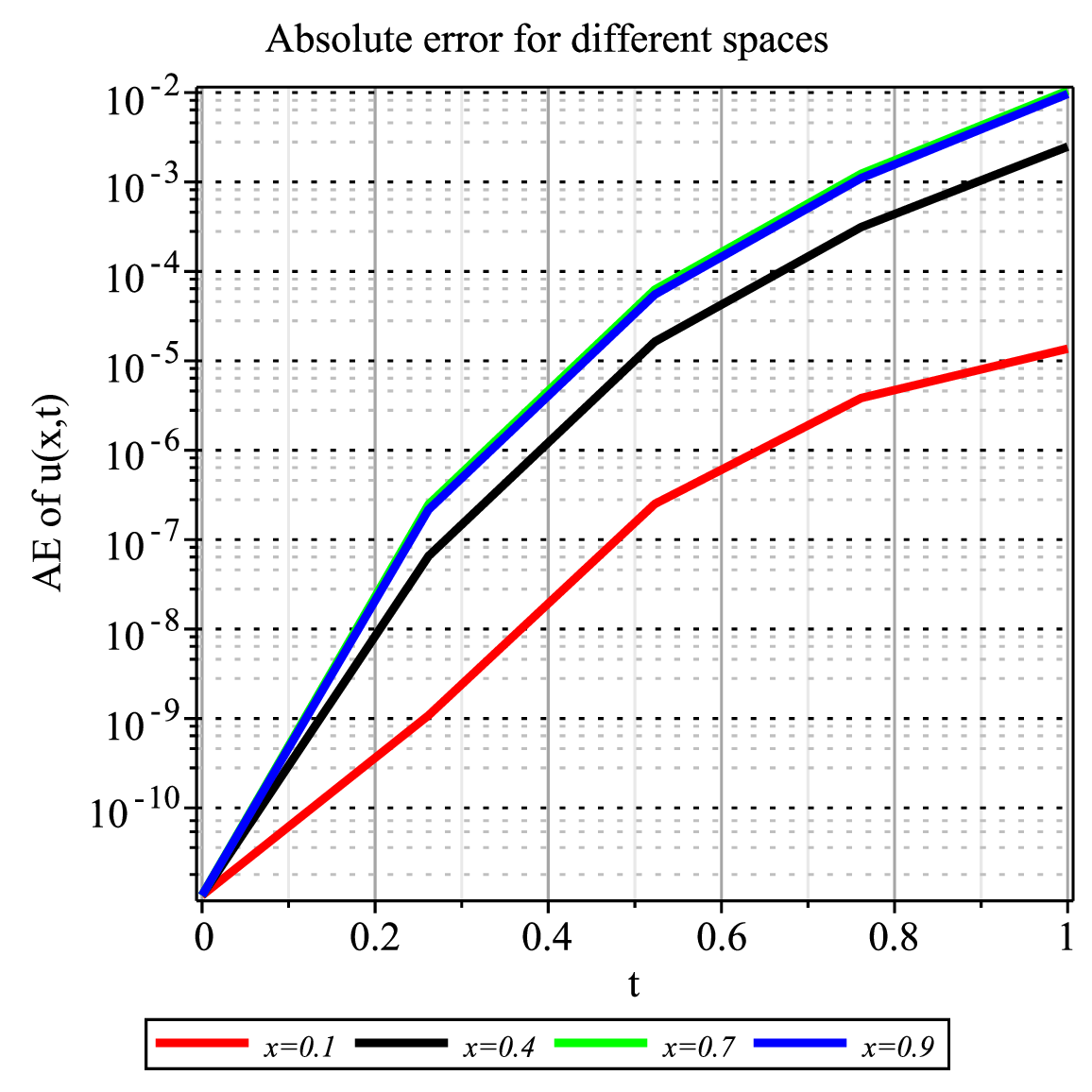}\hspace{0.1em}
    \caption{  {Comparison plots of Exact, Classical, and Modified LADM solutions, along with their corresponding absolute errors for Problem \ref{problem4}.}}\label{fig11}
  \end{figure}
\begin{figure}[H]
    \centering
     \includegraphics[width=5.45cm]{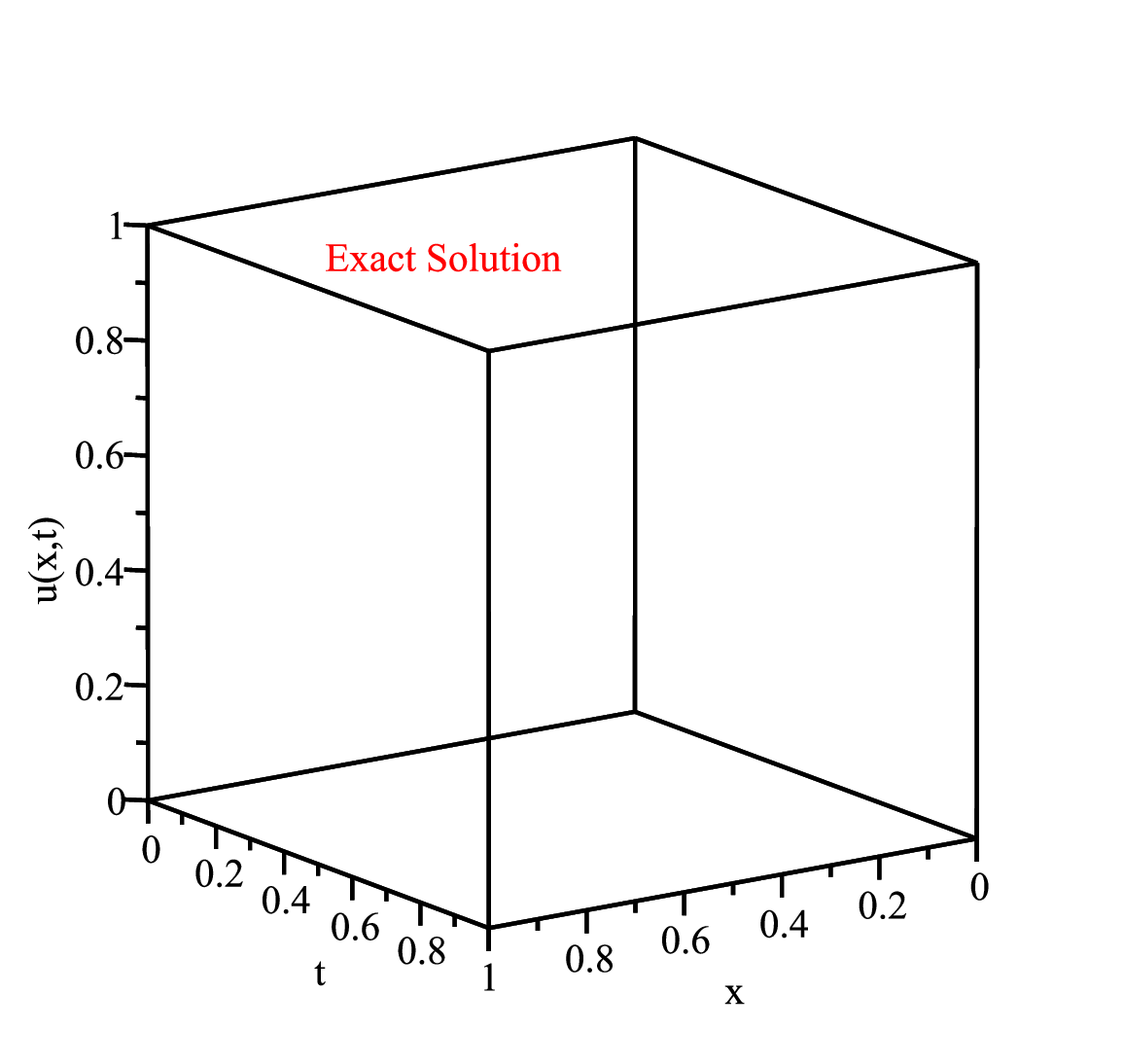}\hspace{0.00em}%
     \includegraphics[width=5.45cm]{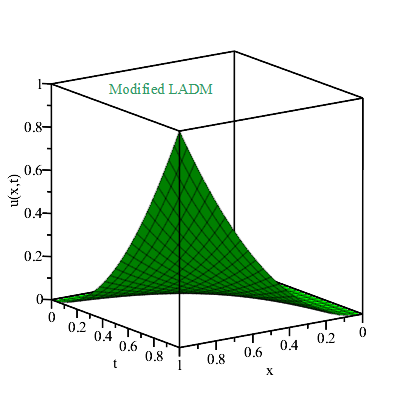}\hspace{0.00em}
     \includegraphics[width=5.45cm]{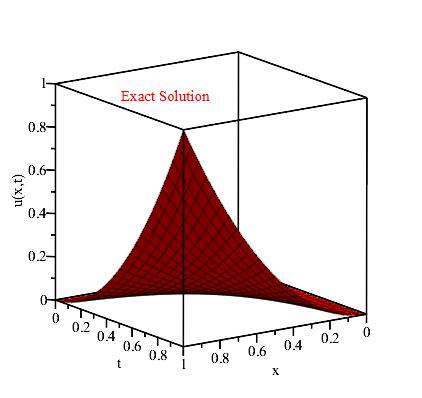}\hspace{0.000em}
    \caption{3D comparison plot presenting the Exact, Classical, and Modified LADM solutions for  problem \ref{problem4}.}\label{fig12}
  \end{figure}
\end{prb}

\begin{prb}\label{nonproblem1} 
Consider the following non-linear initial and boundary value problem
\begin{equation}
{\frac {\partial^\alpha }{\partial t^\alpha }}u \left( x,t \right) +u \left( x,t
 \right) {\frac {\partial }{\partial x}}u \left( x,t \right) -u
 \left( x,t \right) {\frac {\partial ^{2}}{\partial {x}^{2}}}u \left(
x,t \right) =
 h({{{x}}}, {{{t}}})+4\,u(x,t){\pi}^{2}
{t}^{2}\sin \left( 2\,\pi\,x \right),
\end{equation}
with the IBCs condition as follows:
\begin{equation*}
\begin{split}
u(x, 0) &= 0,\\
u(0, t) &=  0, \\
u(1, t) &=0,
\end{split}
\end{equation*} 
where $h(x,t)$ is given by
\begin{eqnarray} \nonumber
h({{{x}}}, {{{t}}})=2\,t\sin \left( 2\,\pi\,x \right) +2\,\pi\,{t}^{4}\sin
 \left( 2\,\pi\,x \right) \cos \left( 2\,\pi\,x \right),\\ \nonumber
\textnormal{for} \quad t > 0, \quad 0 \leq x \leq 1, \quad 0 < \alpha \leq 1.\\\nonumber
\end{eqnarray} 
The exact solution  for  \ref{nonproblem1} is given by
\begin{equation*}
u(x, t) = {t}^{2}\sin \left( 2\,\pi\,x \right).
\end{equation*}

\begin{figure}[H]
    \centering
    \includegraphics[width=5.4cm]{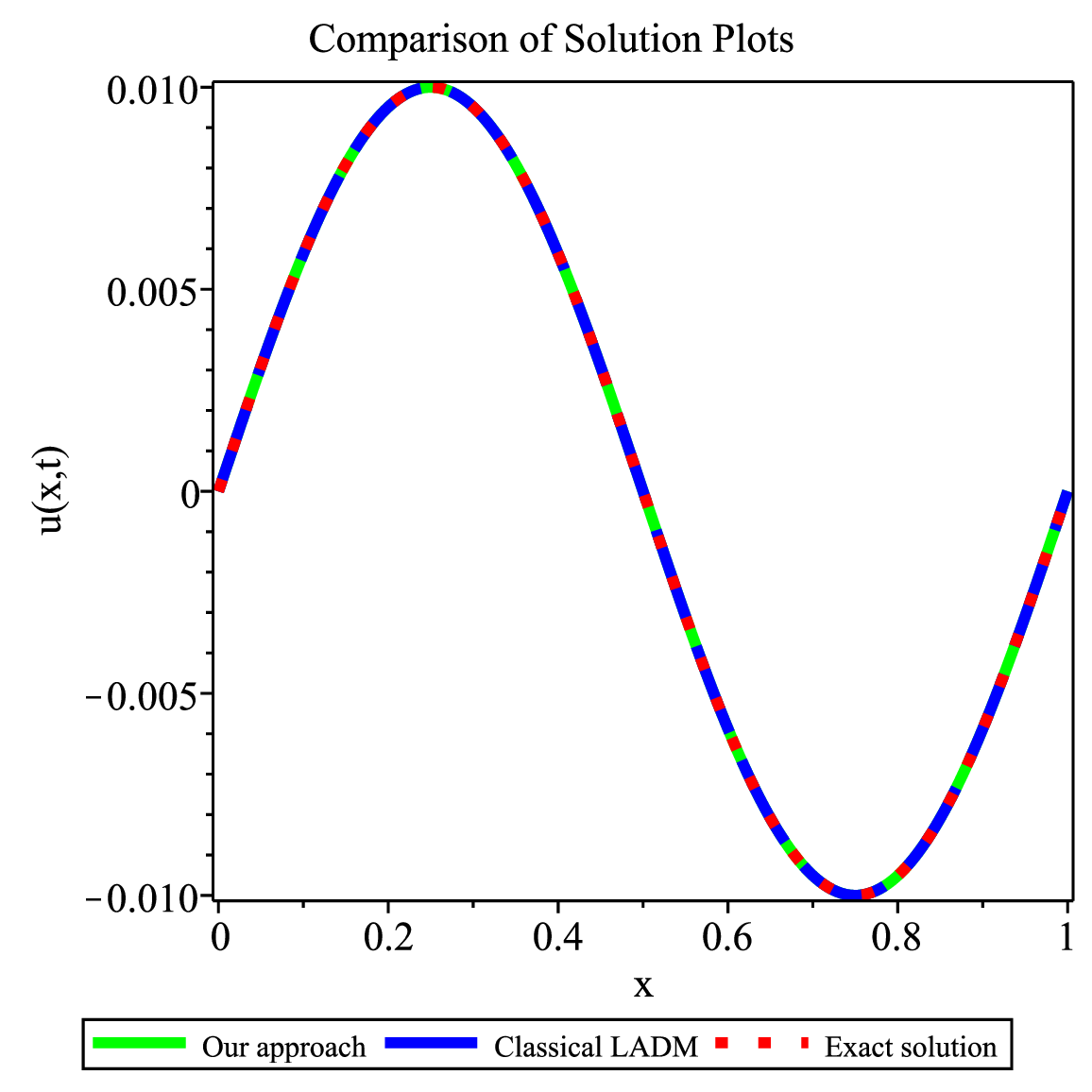}\hspace{0.1em}%
        \includegraphics[width=5.4cm]{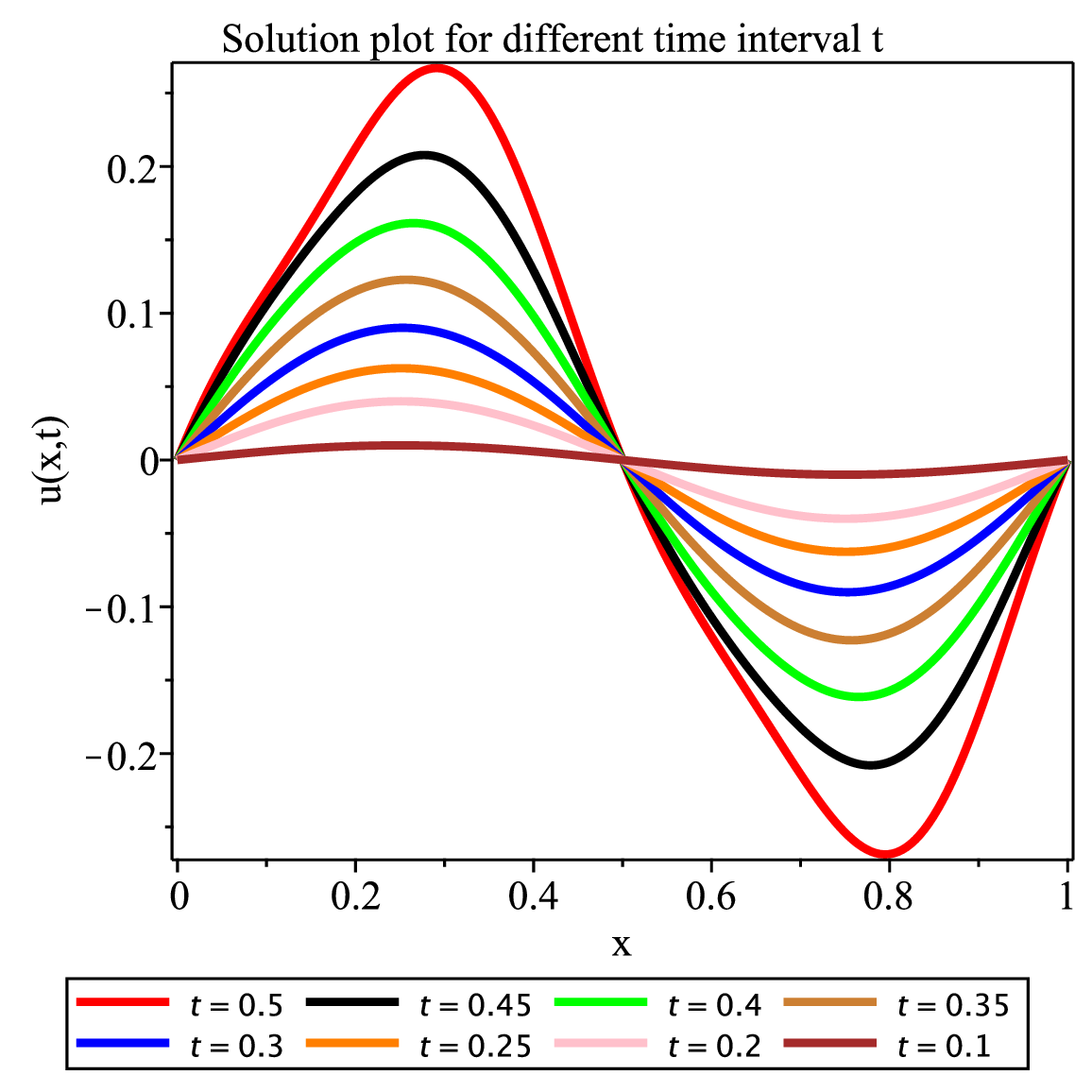}\hspace{0.1em}%
     \includegraphics[width=5.4cm]{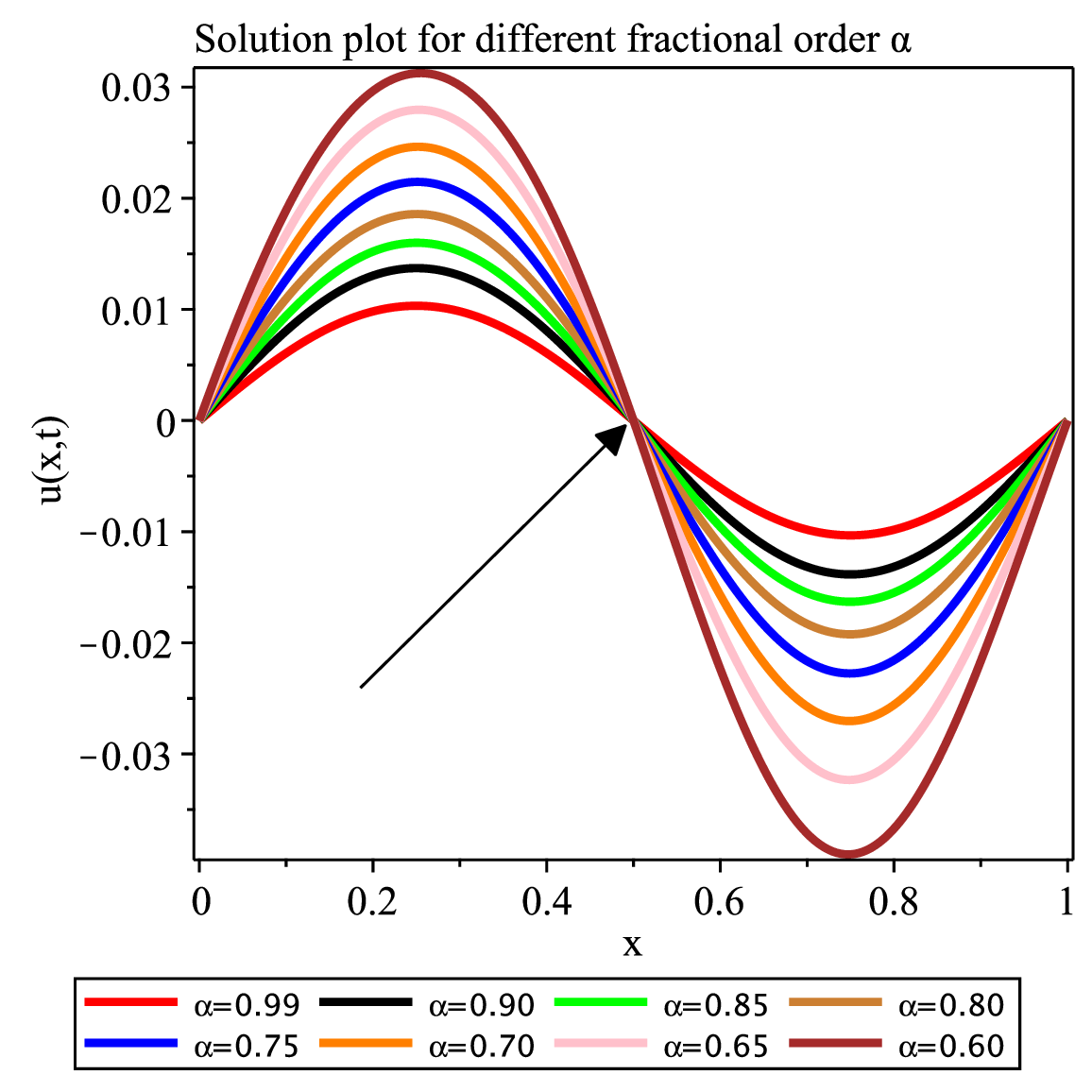}%
    \llap{\raisebox{1.21cm}{%
      \includegraphics[width=1.90cm]{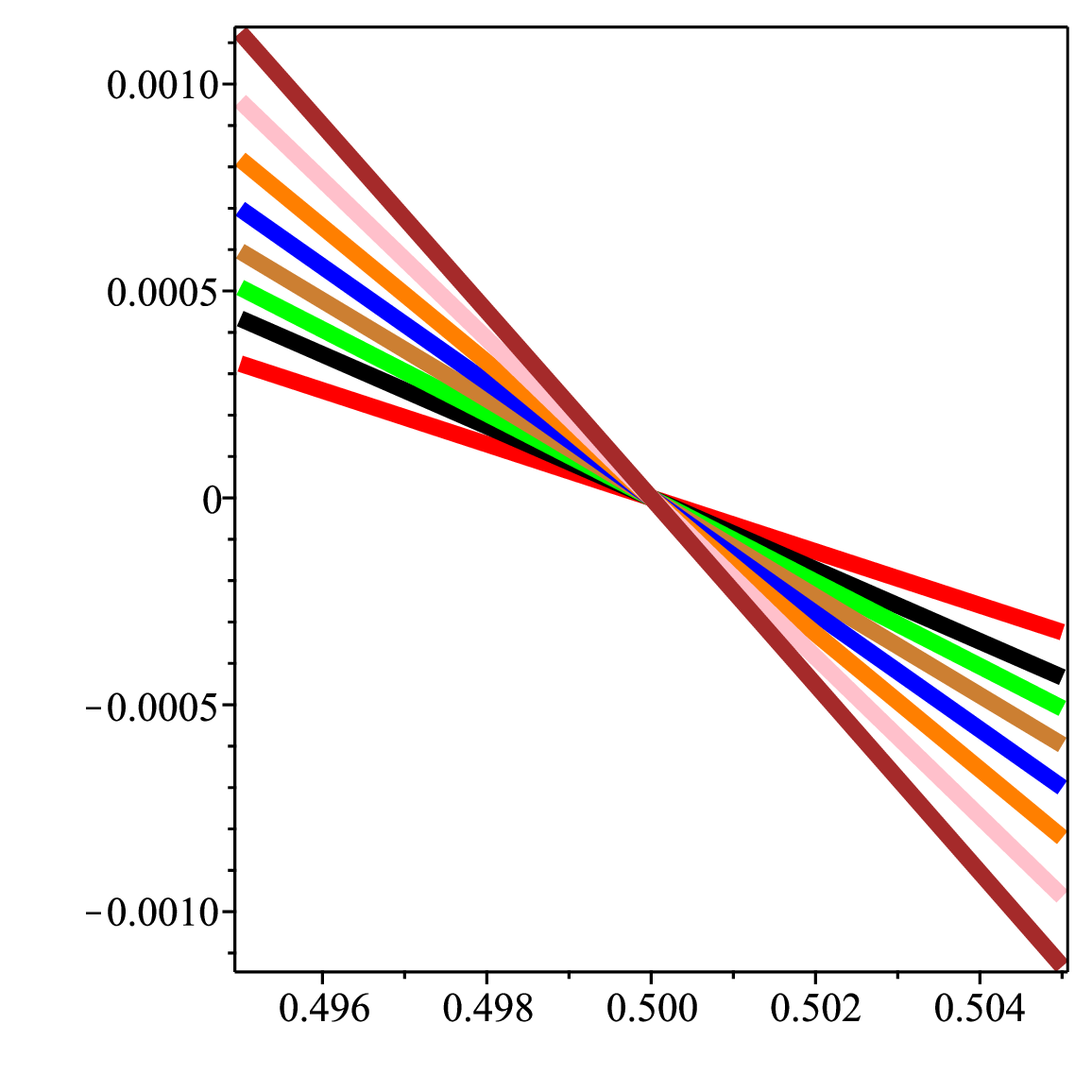}\hspace{6.10em}%
    }}\hspace{0.1em}
    \caption{  {Comparison plots of Exact, Classical, and Modified LADM solutions   at different fractional order $\alpha$ for Problem \ref{nonproblem1}.}}\label{fig13}
  \end{figure}
\begin{figure}[H]
    \centering
     \includegraphics[width=5.45cm]{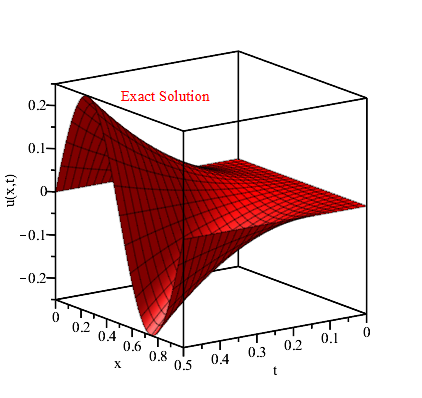}\hspace{0.00em}%
     \includegraphics[width=5.45cm]{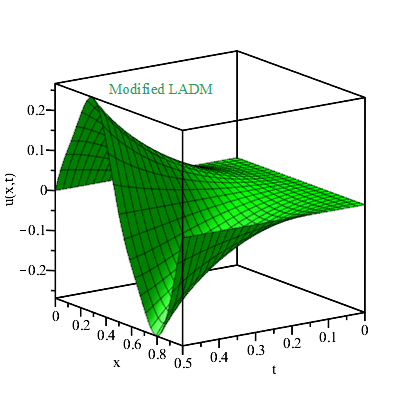}\hspace{0.00em}
     \includegraphics[width=5.45cm]{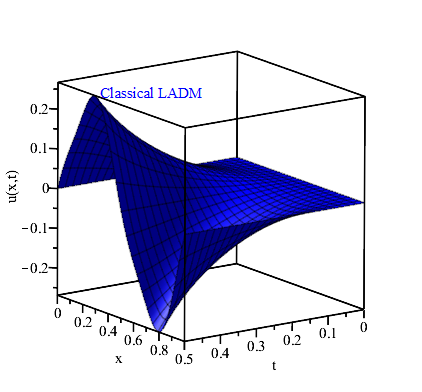}\hspace{0.000em}
    \caption{3D comparison plot presenting the Exact, Classical, and Modified LADM solutions for  problem \ref{nonproblem1}.}\label{fig14}
  \end{figure}
\end{prb}

\section{Result   and Discussion }\label{sec4}

The results obtained from applying the Modified Laplace   Decomposition Method (MLDM) to the various test problems demonstrate its effectiveness and reliability in solving fractional partial differential equations (FPDEs). The convergence analysis shows that the proposed technique effectively approximates the solutions of the considered FPDEs. In Section~\ref{sec3}, we conducted several numerical experiments with varying orders of the fractional derivatives and initial conditions, confirming that the MLDM exhibits a high rate of convergence toward the exact solutions. We provide further analysis on each test as follows.

In Test \ref{problem1}, we investigated a one-dimensional fractional partial differential equation.
The results, illustrated in Figure \ref{fig1}, show that the MLDM solutions align closely with the exact solution, particularly in the region of interest. The absolute error plots indicate that the MLDM significantly reduces the error compared to the classical LADM, confirming the improved accuracy of our modified approach. Furthermore, the 3D comparison plot in Figure \ref{fig2} highlights the behavior of the solutions across the spatial domain, showcasing the effectiveness of the MLDM. In Test \ref{problem2}, we extended our analysis to a two-dimensional case with non-homogeneous boundary conditions. The comparison plots in Figure \ref{fig3} illustrate that both the MLDM and classical LADM methods yield similar results; however, the MLDM consistently exhibits better accuracy at various fractional orders $\alpha$. The 3D representation in Figure \ref{fig4} provides a clear visualization of the solution's structure, reinforcing the advantages of the MLDM in handling more complex problems. Test \ref{problem32} involved an initial boundary value problem with a known analytical solution. The results presented in Figure \ref{fig5} demonstrate that the MLDM effectively captures the exact solution while outperforming the classical method at different fractional orders $\alpha$. The corresponding 3D plot in Figure \ref{fig6} emphasizes the robustness of the MLDM, illustrating how it maintains accuracy across the solution domain. In Test \ref{newp3}, we tackled another initial boundary value problem. The results in Figure \ref{fig7} indicate that the MLDM solutions provide a better fit to the exact solution compared to the classical LADM, particularly near the boundaries. The comparison plots also showcase the corresponding absolute errors, highlighting the MLDM's capability to handle such complexities. The 3D plot in Figure \ref{fig8} illustrates the behavior of the solution, affirming the MLDM's effectiveness in solving this problem.

Test \ref{problem3} examined an initial boundary value problem with zero initial conditions. The results, shown in Figure \ref{fig9}, indicate that the MLDM converges rapidly to the exact solution, with significantly lower errors than the classical method. The 3D comparison in Figure \ref{fig10} further corroborates these findings, showcasing the smoothness of the solution over the domain. In Test \ref{problem4}, we considered a non-linear initial boundary value problem. The results in Figure \ref{fig11} demonstrate that the MLDM achieves high accuracy, particularly as the number of terms in the series solution increases. The corresponding absolute error analysis indicates that the MLDM minimizes discrepancies effectively, as shown in the 3D plot in Figure \ref{fig12}.

Lastly, in Test \ref{nonproblem1}, we analyzed a more complex non-linear initial boundary value problem. The results presented in Figure \ref{fig13} clearly demonstrate that the MLDM maintains accuracy across different fractional orders $\alpha$. The 3D representation in Figure \ref{fig14} highlights the dynamic behavior of the solution, confirming the versatility of the MLDM in addressing non-linearities.

Overall, the numerical experiments confirm that the Modified Laplace Adomian Decomposition Method is a powerful tool for solving initial-boundary value problems associated with fractional dynamics. The results indicate that this approach not only converges rapidly to the exact solutions but also provides high accuracy with minimal computational effort. Future research will focus on extending this method to more complex FPDEs and exploring its applicability in various scientific domains. Our findings contribute to the ongoing efforts to enhance the modeling capabilities of fractional calculus in addressing real-world problems.

\section{Conclusion}\label{conclusion}
In our present work, a modified technique based on the Laplace transformation and decomposition method is implemented to solve initial boundary value problems of fractional partial differential equations. The generalized scheme of the proposed technique is successfully constructed with the help of Caputo operator. The derived scheme is then used to solve several problems of different kinds. The fractional and integral order results are calculated for each problem, which provide useful information about the actual dynamics of the proposed problems. A comparison between the exact solutions, classical LADM, and our approximate solutions is performed using graphs. The graphical representation confirms the validity of the present technique. Our technique utilizes the boundary conditions significantly by employing a double approximation formula in one iteration. The simulation procedure implies a reduced volume of calculations and a higher degree of accuracy for the modified technique. We believe that our newly modified approach can be extended to other linear and non-linear IBVPs, for example Euler equation \cite{HSIAO1999280, HP2000}, Boussinesq Equations \cite{MR1700862, MR2346436}, multi-phase and mixing flows \cite{antontsev1989boundary, MR1700669}, etc. 
\section{Declaration}
 \subsection*{CRediT authorship contribution statement}
Both authors contributed equally to this work and agreed on the final version of it.
\subsection*{Declaration of competing interest}
				The authors declare that they have no known competing financial interests or personal relationships that could have appeared to influence the work reported in this research paper.
				\subsection*{Data availability}
This paper contains no hidden data. Q. Khan developed and implemented the codes using Maple Version 2024 on a GPU system powered by an Intel(R) Core(TM) i9-8950HK  at the STEM lab of the Education University of Hong Kong. Data for different fractional orders, as well as specific values of ${{{t}}}$ and {{{$x$}}}, will be provided upon reasonable request.
				\subsection*{Funding}
This work is partially supported by Hong Kong General Research Fund (GRF) grant project number 18300821, 18300622 and 18300424, and the EdUHK Research Incentive Award project titled ``Analytic and numerical aspects of partial differential equations''.
				\subsection*{Acknowledgment }
We are very grateful to our MIT lab mates and the technicians at The Education University of Hong Kong for providing all the necessary facilities for this research work.
\bibliographystyle{amsalpha}
\bibliography{mybib}
\end{document}